%% file: resolvante2.tex
\documentclass[a4paper]{amsart}

\usepackage{amssymb,amsmath,graphics,verbatim}
\usepackage{latexsym}
\usepackage{eucal}

\makeatletter
\@addtoreset{equation}{section}\makeatother

\newtheorem{theo}{Theorem}[section]
\newtheorem{lem}[theo]{Lemma}
\newtheorem{prop}[theo]{Proposition}

\theoremstyle{remark} \newtheorem{remark}[theo]{Remark}

\newcommand{\mc}{\mathcal}
\newcommand{\rr}{\mathbb{R}}
\newcommand{\nn}{\mathbb{N}}
\newcommand{\cc}{\mathbb{C}}

\newcommand{\zz}{\mathbb{Z}}

\newcommand{\la}{\lambda}
\newcommand{\eps}{\epsilon}

\newcommand{\pl}{\partial}
\newcommand{\x}{\times}

\newcommand{\til}{\widetilde}

\newcommand{\cjd}{\rangle}
\newcommand{\cjg}{\langle}

\newcommand{\demi}{\frac{1}{2}}
\newcommand{\ndemi}{\frac{n}{2}}

\newcommand{\zf}{\textrm{zf}}
\newcommand{\bfo}{\textrm{bf}_0}
\newcommand{\rbo}{\textrm{rb}_0}
\newcommand{\lbo}{\textrm{lb}_0}
\newcommand{\lb}{\textrm{lb}}
\newcommand{\rb}{\textrm{rb}}
\newcommand{\bfa}{\textrm{bf}}
\newcommand{\bfc}{\textrm{bf}}
\newcommand{\sca}{\textrm{sc}}
\newcommand{\ff}{\textrm{ff}}
\newcommand\Id{\operatorname{Id}}
\def\qed{\hfill$\square$}

\newcommand\RR{\mathbb{R}}

\newcommand\extunion{\overline{\cup}}
\newcommand\Omegab{\tilde \Omega_b^{1/2}}
\newcommand\MMksc{M^2_{k, \sca}}

\newcommand\Omegabht{\tilde \Omega_b^{1/2}}

\newcommand\Mbar{M}

\begin{document}
\title[Resolvent at low energy and Riesz transform]{Resolvent at low energy and Riesz transform 
for Schr\"odinger operators on asymptotically conic manifolds. II.}
\author{Colin Guillarmou}
\address{Laboratoire J. Dieudonn\'e\\
Universit\'e de Nice\\ Parc Valrose\\
06100 Nice\\ FRANCE}     
\email{cguillar@math.unice.fr}
\author{Andrew Hassell}
\address{Department of Mathematics, Australian National University \\ Canberra ACT 0200 \\ AUSTRALIA}
\email{hassell@maths.anu.edu.au}

\subjclass[2000]{}
%

\begin{abstract} Let $M^\circ$ be a complete noncompact manifold and $g$ an asymptotically conic metric on $M^\circ$, in the sense that  $M^\circ$ compactifies to a manifold with boundary $M$ in such a way that $g$ becomes a scattering metric on $M$. A special case of particular interest is  that of asymptotically Euclidean manifolds, where $\partial M = S^{n-1}$ and the induced metric at infinity is equal to the standard metric. We study the resolvent kernel $(P + k^2)^{-1}$  and Riesz transform of the operator $P = \Delta_g + V$, where $\Delta_g$  is the  positive Laplacian associated to $g$ and $V$ is a real potential function $V$ that is smooth on $M$ and vanishes to some finite order at the boundary. 

In the first paper in this series we made the assumption that  $n \geq 3$ and that $P$ has neither zero modes nor a zero-resonance and showed (i) that the resolvent kernel is conormal to the lifted diagonal and polyhomogeneous at the boundary  on a blown up version of $M^2 \times [0, k_0]$, and
(ii)  the Riesz transform of $P$ is bounded on $L^p(M^\circ)$ for $1 < p < n$, and that this range is optimal unless $V \equiv 0$ and $M^\circ$ has only one end. 

In the present paper, we perform a similar analysis assuming again $n \geq 3$ but allowing zero modes and zero-resonances. We show that once again that (unless $n=4$ and there is a zero-resonance) the resolvent kernel is polyhomogeneous on the same space and compute its leading asymptotics. This generalizes results of Jensen-Kato and Murata to the variable coefficient setting. We also find the precise range of $p$ for which the Riesz transform (suitably defined) of $P$ is bounded on $L^p(M)$ when zero modes (but not resonances, which make the Riesz transform undefined) are present. Generically the Riesz transform is bounded for $p$ precisely in the range $(n/(n-2), n/3)$, with a bigger range possible if the zero modes have extra decay at infinity. 
\end{abstract}

\maketitle
\section{Introduction}

This is the second in a series of papers on the analysis of low energy  asymptotics of the resolvent for Laplace-type operators on asymptotically conic spaces. Our setting is a complete noncompact Riemannian manifold $(M^\circ, g)$ which is asymptotically conic in the sense that $M^\circ$ is the interior of a manifold with boundary $M$ such that in a collar neighbourhood 
$[0,\eps)_x\x\pl M$ near $\pl M$, 
\begin{equation}\label{metricconic}
g=\frac{dx^2}{x^4}+\frac{h(x)}{x^2}
\end{equation}
where $x$ is a smooth function that defines the boundary $\pl M$ (i.e. $\pl M=\{x=0\}$ and $dx$ does not vanish on $\pl M$) and $h(x)$ is a smooth family of metrics on $\partial M$. 
We let $V$ be a real potential such that there exists an integer $l\geq 3$ with 
\begin{equation}\label{hyp2}
V\in C^{\infty}(M), \quad V=O(x^{l}) \textrm{ as }x\to 0, 
\end{equation}
and consider the Schr\"odinger operator   
\[P=\Delta_g+V\]
where $\Delta_g$ is the positive Laplacian with respect to $g$. Our main interest is the behaviour of the resolvent kernel $(P + k^2)^{-1}$ as $k \downarrow 0$, and related operators such as the Riesz transform and the heat kernel of $P$.

The operator $P$
is self-adjoint on $L^2(M,dg)$ and its spectrum is $\sigma(P)=[0,\infty)\cup \sigma_{\rm pp}(P)$
where $\sigma_{\rm pp}(P)=\{-k^2_1\geq \dots \geq-k^2_N\}$ is a set of negative eigenvalues (by convention $k_i>0$). The point $0$ can
also be an $L^2$-eigenvalue, but the positive spectrum is absolutely continuous. 

The resolvent $R(k)=(P+k^2)^{-1}$ is well defined as a bounded operator on $L^2(M,dg)$
for $k \in (0, k_1)$ but fails to be bounded or defined at $k=0$. 
In the first paper in this series, \cite{GH}, which we refer to as Part I, we assumed that $n \geq 3$, that $l = 3$ in \eqref{hyp2} and that  $P$ had neither zero modes (i.e. $0$ is not an $L^2$ eigenvalue of $P$) nor a zero-resonance and showed that for small $k_0 < k_1$,  the  distributional kernel of $R(k)$ 
is well-behaved on a certain manifold with corners, denoted $\MMksc$, that is a blown-up version of $X=[0,k_0]_k\x M\x M$. 
This space, illustrated in Figure ~\ref{mmksc}, is a compact manifold with corners up to codimension 3, and has 8 boundary hypersurfaces. Each boundary hypersurface can be considered an `asymptotic regime'; for example, one asymptotic regime is when $k$ goes to zero and simultaneously  the left variable $z$, $z \in M$, tends to infinity (i.e. $x(z)\to 0$) while $z'$ is held fixed, and this corresponds to the boundary hypersurface $\lbo$ in the figure. This face $\lbo$ may be considered as a bundle  of directions of convergence in the corner $\{x=0,k=0\}$, the fibers of which are identified with
$[-1,1]_\tau$ where $\tau:=(x-k)/(x+k)$, the basis being the corner.   

Using this space we defined a calculus of pseudodifferential operators, 
denoted $\Psi_k^{m, (a_{\bfo}, a_{\zf}, a_{\sca}); \mc{E}}(M; \Omegabht)$ which are polyhomogeneous conormal both at the diagonal submanifold and at each boundary hypersurface of $\MMksc$. The index $m$ denotes the conormal order at the diagonal (that is, the order as a pseudodifferential operator), while $(a_{\bfo}, a_{\zf}, a_{\sca})$ and $\mc{E}$ specify the type of polyhomogeneous expansion allowed at the various boundary hypersurfaces of $\MMksc$; see Section 2.3 of Part I. The first main result of Part I was that the resolvent kernel lies in this calculus of operators; in this sense we found the complete asymptotic expansion of the resolvent kernel, in every asymptotic regime,  as $k \to 0$. We also computed the leading order behaviour of the resolvent at each boundary hypersurface. 

This precise analysis of the resolvent at low energy allowed us to determine the exact range of $p$ for which the Riesz transform of $P$ is bounded on $L^p$. If $P$ is nonnegative, then the Riesz transform $T$ is defined by $T = d P^{-1/2}$; in general, to make sense of this, one needs to project off the nonpositive part of the spectrum (i.e. that corresponding to eigenvalues) before taking the $-1/2$ power of $P$. 
In Part I, we showed 

\begin{theo}
Let $n \geq 3$, and let $P = \Delta_g$ be the Laplacian on an asymptotically conic manifold $(M,g)$. First we assume that  $M$ has one end.
Let $\Delta_{\pl\Mbar}$ be the Laplacian on the boundary of $\Mbar$ for the 
metric $h(0)$ given by (\ref{metricconic}),  let $\lambda_1$ be its first non-zero eigenvalue, and let $\nu_1 = \sqrt{((n-2)/2)^2 + \lambda_1}$. 
If $\nu_1<n/2$, then the Riesz transform $T$ satisfies 
\[ T\in \mc{L}(L^p(M),L^p(M;T^*M)) \iff 1<p<n \Big(\frac{n}{2} - \nu_1\Big)^{-1} \] 
while if $\nu_1\geq n/2$, then
 \[ T\in \mc{L}(L^p(M),L^p(M;T^*M)),\quad \textrm{ for all }\quad  1<p<\infty.\] 

Next let $P = \Delta_g + V$ with $V$ as above, and suppose that either $M$ has more than one end, or that $V$ is not identically zero. In either case we assume that $P$ has no zero modes or zero-resonances.  Then
\[ T\in \mc{L}(L^p(M),L^p(M;T^*M)) \iff 1<p<n.\]
\end{theo}

In this second paper, we allow $P$ to have zero modes or zero-resonances. Here a zero mode is an eigenvector $\psi\in L^2(M)$ such that $P\psi=0$, and a zero-resonance is an function $\psi \notin L^2(M)$ with $P\psi = 0$, with $\psi$ tending to zero at infinity; zero-resonances can only occur in dimensions 3 and 4. We carry out a similar parametrix construction for the resolvent as in Part I \cite{GH}, but with additional complications. The main problem is that the leading terms in the parametrix are more singular than they were when no zero modes are present: for example, the leading term is at order $-2$ instead of $0$ at $\zf$, and may be at order $n/2 - 4$ instead of $n/2 - 2$ at $\rbo$. This means we must define several `models' (i.e. terms in the Taylor series) at each boundary hypersurface, rather than just one, and this 
in turn means the compatibility conditions between models on adjacent  boundary hypersurfaces are much more involved. 
Nevertheless, we carry out the program and show that the resolvent lies in our calculus of operators, except in exceptional cases. The result depends on the nature of the $L^2$ kernel of $P$. Define 
\begin{equation}
m' = \min(2, m), \quad m = \sup \{ a \geq 0 \mid \ker P \subset x^{n-2+ a} L^\infty(M) \}.
\label{m'-defn}\end{equation}
We show

\begin{theo} Let $(M,g)$ be asymptotically conic of dimension $n \geq3$. Let $P = \Delta_g + V$ as above and assume that $P$ has nontrivial $L^2$ kernel. Also, if $n=3,4$ or $5$ we assume that 
\begin{equation}
\text{$P$ has no zero-resonance, and $m'>(5-n)/2$ if $n=3,4,5$.}
\label{m'-cond}\end{equation} 
Then $R(k) = (P + k^2)^{-1} \in \Psi_k^{-2; (-2, 0, 0), \mc{R}}(M; \Omegabht)$ is in the calculus and has leading behaviour at order $-2$ at $\zf$, $-2$ at $\bfo$, $n/2 - 4 + m'$ at $\lbo$ and $\rbo$, $0$ at $\sca$ and is rapidly vanishing at $\lb, \rb, \bfc$. 
\end{theo}

This is a simplified version of Theorem~\ref{res-conic-1} which contains additional information. In the asymptotically Euclidean, rather than conic, setting the corresponding result is presented in Theorem~\ref{thm:nullspace}.

The complete analysis in dimensions 3 and 4 is more involved because of the possibility of resonant states and due to the relatively slow decrease at infinity of zero modes. We have not attempted a complete treatment of all cases, but deal fully with the asymptotically Euclidean case in dimension $3$ in Section~\ref{dim3}, which allows a direct comparison with the work of Jensen-Kato \cite{JK}. We also show polyhomogeneity in this case; see Theorem \ref{resolventn=3}. 
The asymptotically \emph{conic} case in dimension $3$ is sketched in Section \ref{sec6} assuming that there is either a unique resonant state 
or a unique zero mode with rather weak decrease at infinity. This assumption is made to avoid excessive technicalities and suffices to  show a variety of possible asymptotic behaviours of the resolvent as $k \to 0$, depending on the rate of decrease of the resonant state (resp. zero mode) at infinity.  
It turns out that in dimension 3, when there is a resonant state $\psi \in x^{3/2}L^\infty(M)\setminus x^{3/2+\eps}L^\infty(M)$, thus failing only logarithmically to be a zero mode, the resolvent is no longer  polyhomogeneous on $M^2_{k,\sca}$. Actually we do not 
give details in this case since the same phenomenon appears in the dimension $4$ asymptotically Euclidean case 
for a resonant state decreasing like $x^2$ (i.e. $r^{-2}$) at infinity, and we have thus chosen to explain this phenomenon in that 
more familiar geometric setting --- see Theorem \ref{n4resonance}. 
It leads to the surprising (although known since at least Murata's work \cite{Mu}) result that the generalized projection onto the resonance occurs at order $k^{-2} (\log k)^{-1}$. The reciprocal of the logarithm shows that the resolvent is not polyhomogeneous in this case; rather it has  the form of $k^{-2}$ times a  convergent power series expansion in $1/\log k$ as $k \to 0$. 

It is worthwhile to emphasize that the analysis of the expansion at $k=0$ of the resolvent for perturbations of conic manifolds has 
been worked out by X-P. Wang \cite{XPW}, but our result is stronger in the sense that we study the 
Schwartz kernel at all asymptotic regimes and not only at $k=0$ with $z,z' \in M^\circ$ fixed. In particular our result about the resolvent allows  us
to obtain precise estimates on the heat kernel as $t\to \infty$ (as shown for instance in \cite{CCH}) and Riesz transform analysis.\\ 
 
We now turn to a generalization of the Riesz transform results to this setting. 
Let $\Pi_>$ denote the spectral projection for the operator $P$ onto the positive spectrum $(0, \infty)$, and $P_> = P \circ \Pi_>$. 
The operator $P_>^{-1/2}$ is obtained from the resolvent using the formula
\begin{equation}\label{p-demi}
 (P_>)^{-\demi} = \Big( \frac{2}{\pi} \int_0^\infty (P + k^2)^{-1}\circ \Pi_>  \, dk \Big) .
\end{equation}
If $P$ has zero modes, then there will be a term  $k^{-2} \Pi_0$ present in the resolvent, which obviously makes the integral diverge. However, composing with $\Pi_>$ kills the $k^{-2}$ singular term and then there is a chance that the integral becomes convergent,  allowing  the Riesz transform to be defined in this setting too. It turns out that this happens in dimensions $n \geq 6$, but not necessarily in lower dimensions where terms of order $k^{-\alpha}$ with $\alpha\geq 1$ may appear in the expansion of $R(k)$ at $k=0$. In the present paper, we find the optimal range of $p$ for which the Riesz transform is bounded on $L^p$; this range depends on $n$ and the nature of the null space through the number $m'$ defined in \eqref{m'-defn}. There is a large literature about this question
for the Laplacian on manifolds (i.e. when $V=0$) and a few results for the case $V\geq 0$; we refer
the reader to the Introduction of Part I \cite{GH} for a few references and known results.

Since we need to compute explicitly the leading asymptotic coefficient of the Schwartz kernel of the resolvent $R(k;z,z')$
as $z'\to 0$ to analyze Riesz transform, we
will make an additional assumption which is not serious but more a technical device to keep the paper
of a reasonable clarity: we shall suppose that 
\begin{equation}\label{assump1}
h(x)-h(0)=O(x^3), \quad V=O(x^5) \textrm{ as }x\to 0.
\end{equation}
\begin{theo}\label{Riesz}
Let $(M,g)$ be an asymptotically conic manifold of dimension $n \geq 3$ and $V$ be a potential satisfying \eqref{assump1}. With $P:=\Delta_g+V$,
assume that $\ker_{L^2} P \not=0$ and \eqref{m'-cond}, where $m'$ be defined by \eqref{m'-defn}. Then  
\[ 
T\in \mc{L}(L^p(M),L^p(M;T^*M))\iff \frac{n}{n-2+m'}<p<\frac{n}{3-m'}.
\]
\end{theo}

\begin{remark}
In dimension $3$ and when $1<m\leq 3/2$, we see that $T$ is bounded in the range $ (3/(1+m),3/(3-m))\subset (1,2)$ 
and when $m\to 1$, this interval reduces to $p=3/2$. It is a bit surprising that 
it does not contain the middle range $2$, but this regime of $m$ is somehow special: we  
prove directly (without using resolvent) in Section \ref{unboundedness} that, in this case, there is a sequence $f_R\in \Pi_> (L^2(M))$
such that $||df_R||^2 \geq \alpha(R) \cjg Pf_R,f_R\cjd $ with $\alpha(R)\to \infty$ when $R\to \infty$. 
The condition $m\leq 1$ also corresponds to the case where $R(k)\Pi_>$ becomes non-integrable at $k=0$ since
its expansion at $k=0$ contains powers of order lower or equal to $-1$.  
In dimension $4$ (resp. $5$), the same phenomenon holds when $1/2<m<1$ (resp. $m=0$), then the limit as 
$m\to 1/2$ (resp. $m\to 0$) of the interval of boundedness reduces
to $p=8/5$ (resp. $5/3$).    
\end{remark}

\begin{remark}
These cases might be considered somewhat artificial, geometrically speaking, since the resonances or eigenvalues
come from the potential perturbation $V$. However, it is very likely that a similar phenomenon occurs
for the Laplacian acting on forms on asymptotically conic and Euclidean manifolds, due to $L^2$ harmonic forms and zero-resonance forms \cite{GC}. We plan to analyze this in a future publication.   
\end{remark}

\section{The blown-up space $\MMksc$ and the calculus $\Psi_k^{m, \mc{E}}(M; \Omegabht)$.}

\subsection{The space $\MMksc$.}

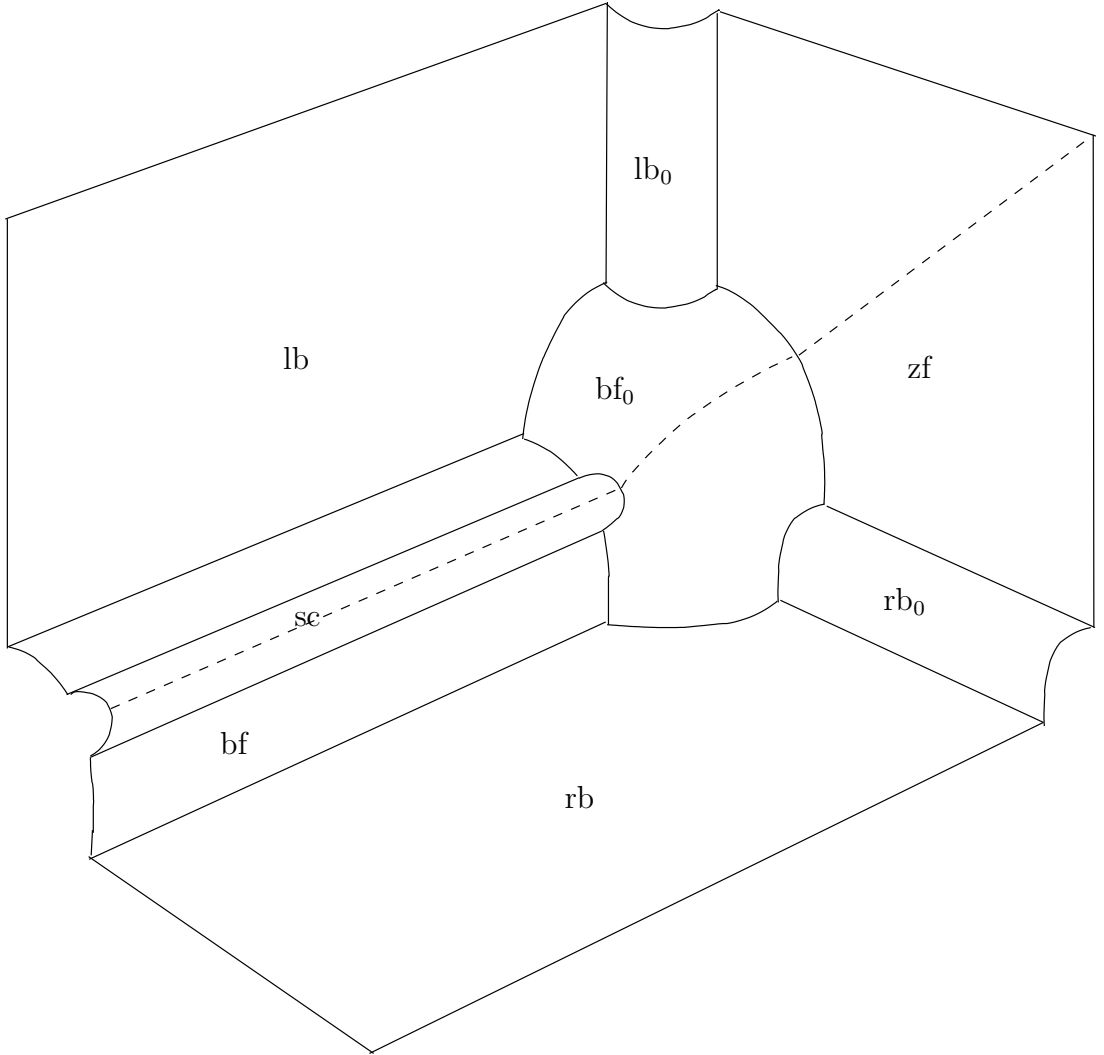
\begin{figure}[ht!]\label{mmksc}
\begin{center}
\input{riesz2.pstex_t}
\caption{The blow-up manifold $M^2_{k,\sca}$.}
\end{center}
\end{figure}

The resolvent kernel $(P + k^2)^{-1}$ is a distribution on the space $M^\circ \times M^\circ \times (0, k_0]$. We work on a compactification of this space tailored to the geometric and analytical properties of the operator $P + k^2$. We begin with the obvious compactification $M \times M \times [0, k_0]$ and perform several blowups. The reason for the blowups can be explained in terms of the off-diagonal behaviour of the resolvent kernel. Roughly speaking at energy $-k^2$, the resolvent kernel should have exponential decay at the scale $k^{-1}$ away from the diagonal. Thus for fixed $k > 0$ the kernel is essentially supported close to the diagonal, but for $k =0$ we get only weak (polynomial) decay. To effectively decouple these two effects we blow up boundary submanifolds of $M^2 \times [0, k_0]$ which separates the exponential behaviour for $k > 0$ from the polynomial behaviour at $k=0$. Let us denote the boundary hypersurfaces of this space $\zf = M^2 \times \{ 0 \}$, $\rb = M \times \partial M \times [0, k_0]$, and $\lb = \partial M \times M \times [0, k_0]$. Then  we blow up the submanifold $(\partial M)^2 \times \{ 0 \}$, followed by the lift to this space of $(\partial M)^2 \times [0, k_0]$, $M \times \partial M \times \{ 0 \}$, $\partial M \times M \times \{ 0 \}$, to produce a space we call $M^2_{k, b}$. The new boundary hypersurfaces so created are denoted $\bfo$, $\bfc$, $\rbo$ and $\lbo$, respectively. 
Finally we blow up the boundary of the diagonal intersected with $\bfc$ to create a new boundary hypersurface $\sca$, to produce the final space $\MMksc$. See Figure~\ref{mmksc}. The rigorous description of this space is done in Subsection 2.2 of Part I \cite{GH}. 

This space has eight boundary hypersurfaces, and each one can be thought of as a geometric realization of a distinct `asymptotic regime'. 
Let us denote points in the left copy of $M$ by $z$ and points in the right copy of $M$ by $z'$. For $z, z'$ in a collar neighbourhood of $\partial M$ we write $z = (x,y)$ with $y$ a local coordinate on $\partial M$. 
Then if $k$ tends to $0$ with $z, z'$ fixed we arrive at $\zf$. If $k \to 0$ and $z, z'$ both tend to infinity with $k/x$ and $k/x'$ tending to limiting values, then we arrive at $\bfo$. If $k > 0$ is fixed and $z, z'$ both tend to infinity with $y$ tending to a limit and the distance $d(z,z')$ fixed then we arrive at $\sca$, and so on. 

To construct an accurate parametrix for the resolvent, we need to solve the PDE $(P + k^2) G = \delta$ where $\delta$, the kernel of the identity function, is a delta function supported on the diagonal submanifold. Thus the construction `begins' by specifying the correct
conormal singularity at the diagonal, together with the correct model operators on each of the boundary hypersurfaces that intersects the diagonal, namely $\sca$, $\zf$ and $\bfo$. Let us consider the model problem at $\zf$, that is at $k=0$ and for $z, z' \in M^\circ$, in more detail. 

We define the metric $g_b$ by $g_b = x^2g$; it is then an exact b-metric in the sense of Melrose \cite{APS}, that is an asymptotically cylindrical metric on $M$. We now define the operator $P_b$ by 
\begin{equation*}\label{pvspb}
P=x^{\ndemi+1}\, P_b \, x^{-\ndemi+1}.
\end{equation*} 
Since $P$ is formally self-adjoint with respect to $g$, $P_b$ is formally self-adjoint with respect to $g_b$. Actually, we prefer to think about the relation between $P$ and $P_b$ differently. If we regard the derivatives comprising $P_b$ as annihilating the b-half-density $|dg_b|^{1/2}$ (the Riemannian half-density corresponding to the metric $g_b$), and the derivatives comprising $P$ as annihilating the half-density $|dg|^{1/2}$, then this effectively implements the conjugation by $x^{n/2}$ in the formula above and the relation becomes
\begin{equation}\label{pvspb**}
P=x \, P_b \, x.
\end{equation} 
Explicitly, the formula for $P_b$ is 
\begin{equation}\label{pb-form}
P_b:=-(x\pl_x)^2+\Big(\frac{n-2}{2}\Big)^2+\Delta_{\pl\Mbar}+
W,
\end{equation}
with $W\in x\textrm{Diff}_b(M)$ a lower-order term at $x=0$ (here $\textrm{Diff}_b$ means the set of  differential operators
which are smooth combinations of vector fields tangent to $\pl M$, see Part I \cite{GH}). 
We see directly from this that $P_b$ is elliptic as an b-differential operator. The kernel of its inverse therefore lives naturally on the space $M^2_b = [M^2; (\partial M)^2]$, i.e. $M^2$ with the corner blown up \cite{APS}. The face $\zf$ is canonically diffeomorphic to $M^2_b$, and the leading model for $(P + k^2)^{-1}$ (at least in the absense of zero-modes and zero-resonances) is the inverse of $P_b$ (up to powers of $x$). 

We need to recall some of the main results of 
the b-calculus from \cite{APS}. For more details, 
and definitions of index sets, etc,  see Part I \cite{GH}, 
section 2; here we just quote the three main theorems needed 
in this paper: 

First, we denote $L^2_b(M):=L^2(M;\Omega_b^\demi)$ the square
integrable (b-)half-densities; here $\Omega_b$ is a 
canonical smooth bundle over $M$
of 1-densities defined in Subsection 2.2.2 of Part I \cite{GH}
and trivialized by $|dg_b|$.
The spaces $H_b^j(M)$ is the associated $j$-th Sobolev space, see Section 2.1 of
Part I.

Let $\lambda_0 = 0 < \lambda_1 \leq \lambda_2 \dots$ be the spectrum of the Laplacian on $(\pl\Mbar,h_0)$. We then define 
\begin{equation}\label{npl}
N_{\pl}:=\{ \nu_0, \nu_1, \dots \mid \nu_i = \sqrt{((n-2)/2)^2 + \lambda_i  } \}.\end{equation}
In the case of the canonical sphere $(\pl\Mbar=S^{n-1},h_0=d\theta^2)$ (i.e. for the Euclidean setting) we have 
\[N_{\pl}=\{(n-2)/2+l; l\in\nn_0\}. \]
We shall also use the notations
\begin{equation*}
\begin{gathered}
E_{\nu_i}:=\ker(\Delta_{\pl M}-\la_i) \textrm{ for the conic case }, \\
\quad E_{i}:=\ker(\Delta_{S^{n-1}}-i(i+n-2)) \textrm{ for the Euclidean case}.\end{gathered}\end{equation*}

\begin{theo}[Melrose's Relative Index Theorem, \cite{APS}]\label{rit} The operator $P_b$ is Fredholm as a map from $x^\alpha H^j_b(M)$  to $x^\alpha H^{j-2}_b(M)$ for all $j \geq 2$ and all $\alpha \neq \pm \nu_l$, $l = 0, 1, 2, \dots$. The index of $P_b$ is equal to $0$ for $|\alpha| < (n-2)/2$ and the index jumps by $d_{l}$, the multiplicity of the $\nu_l^2$-eigenspace $E_{\nu_l}$ of $\Delta_{\pl\Mbar}+(n-2)^2/4$, as $\alpha$ crosses the value $\pm\nu_l$ (with $\alpha$ decreasing).  
\end{theo}

\begin{theo}[Regularity of solutions to $P_b u = f$]\label{reg} Suppose that $f$ is polyhomogeneous on $M$ with respect to the index set $E$, that $u \in x^\alpha L^2_b(M)$, and that $P_b u = f$. Let $S$ be the set
$$S = \{ \pm \nu_l \mid l = 0, 1, 2, \dots \}
$$
and for $b \in \RR$ let $\mu(b,j)=\sharp \{b+k;k=0,\dots j\}\cap S$. 
Then $u$ is polyhomogeneous with respect to the index set $E \extunion F$, where $F$ is the index set 
$$
\big\{ \big((\pm\nu_l+j), k\big) \mid  l\in\nn_0, j\in\nn_0, \ \pm \nu_l\geq \alpha, \ 0 \leq k \leq \mu(\pm\nu_l,j) - 1 \big\}.
$$
When $(\pl\Mbar,h_0)=(S^{n-1},d\theta^2)$, this reduces to 
\[\big\{ \big(n/2+l, k\big) \mid  l\in\zz, \ n/2+l\geq \alpha, \ 0 \leq k \leq N_l - 1 \big\}.\]
where $N_l$ is the number of elements of the form $\pm(n/2-1+j), j\in\nn_0$ in the interval $(\alpha,n/2+l]$ .
\end{theo}

Recall from Part I \cite{GH} the definition of $G_{\nu_l} \subset E_{|\nu_l|}$. 
It is a finite dimensional subspace of $C^\infty(\partial M)$. Informally it is set of coefficients of the leading asymptotic as $x\to 0$ of all $u \in \operatorname{null} P$, $u \in x^{\nu_l - \epsilon}L^2(M)$ for all $\eps>0$. 
\begin{prop}\label{prop:complementary} (\cite{APS}, Chapter 6) The subspaces $G_{\nu_l}$ and $G_{-\nu_l}$ of $E_{\nu_l}$ are orthogonal complements with respect to the inner product on $(\pl\Mbar,h_0)$. 
\end{prop}
In order to simplify the exposition and avoid tedious computations, we will make the assumption \eqref{assump1}. 
This assumption implies that if $u$ is solution of $P_bu=O(x^{\nu+2+\eps})$ with asymptotic 
\[u=\sum_{\substack{\nu_i\in N_{\pl M}\\
\nu\leq \nu_i\leq \nu+2}}x^{\nu_i}(\log x)^{k_i}u_{\nu_i}(y)+O(x^{\nu+2+\eps})\] 
for some $\nu\in N_{\pl M}$, 
$\eps>0$, $k_i\in\nn_0$, $u_{\nu_i}\in C^{\infty}(\pl M)$, then $k_i=0$ and  
$u_{\nu_i}\in E_{\nu_i}$ for any $i$ such that $\nu_i\leq\nu+2$. This is straightforward from the indicial equation.

\subsection{Operator calculus}
The purpose of the space $\MMksc$ is to carry the kernel of our parametrix for the resolvent $(P + k^2)^{-1}$; eventually, we will see that the resolvent kernel itself is nice on this space. We define a `calculus' of pseudodifferential operators, denoted $\Psi_k^{m, \mc{E}}(M; \Omegabht)$, defined by their kernels which are half-densities on $\MMksc$. They can be realized as operators on half-densities on  $M^\circ$ depending parametrically on $k$. The space $\Psi_k^{m, \mc{E}}(M; \Omegabht)$
depends on a pseudodifferential order $m$, controlling the diagonal singularity, and an index family $\mc{E}$, i.e an index set for each of the boundary hypersurfaces $\zf, \bfo, \sca, \lbo, \rbo$, specifying the allowed terms in the polyhomogeneous expansion of the kernel at each of these faces (at the remaining boundary hypersurfaces, i.e. $\lb, \rb, \bfc$, the kernels are required to be rapidly decreasing). See Definition 2.8 in Part I \cite{GH} for the precise definition. 

The most important properties of this calculus are 
\begin{itemize}
\item
There is a composition law: given index families $\mc{A}$ and $\mc{B}$ as above, there is an index family $\mc{C}$ given in Part I \cite[Prop 2.10]{GH} such that
$$
 \Psi_k^{m, \mc{A}}(M; \Omegabht)  \circ \Psi_k^{m', \mc{B}}(M; \Omegabht) \subset \Psi_k^{m+m', \mc{C}}(M; \Omegabht).
 $$
 
 \item If $E \in \Psi_k^{m, \mc{E}}(M; \Omegabht)$ is in the calculus, with $m < 0$ and with $\mc{E}_{\textrm{f}} > 0$ for $\textrm{f} = \zf, \bfo, \sca, \lbo, \rbo$, then for large enough $N$, $E^N$ is Hilbert-Schmidt with $\| E^N(k) \|_{HS} \to 0$ as $k \to 0$. In particular, the operator $\Id - E^N(k)$ is invertible for $N$ large enough and $k$ small enough, and the Neumann series $\Id + E(k) + E(k)^2 + \dots$ for the inverse converges in operator norm and in Hilbert-Schmidt norm. 
 
 \item If the index family $\mc{A}$ is nonnegative, then the restriction of $A \in \Psi_k^{m, \mc{A}}(M; \Omegabht)$ to any of the faces $\zf, \bfo, \sca$ is well-defined, denoted $I_{\textrm{f}}(A)$, $\textrm{f} = \zf, \bfo, \sca$, and called the normal operator at $\textrm{f}$. The kernel $I_{\zf}(A)$ is a b-pseudodifferential operator of order $m$, the kernel $I_{\bfo}(A)$ is a pseudodifferential operator acting on half-densities on $\partial M \times (0, \infty)$ and the kernel $I_{\sca}(A)$ is a family, parametrized by $\partial M \times (0, k_0]$,  of convolution pseudodifferential operators acting on half-densities on $\RR^n$. The 
 normal operators respect composition:
 $$
 I_{\textrm{f}}(A) \circ I_{\textrm{f}}(B) = I_{\textrm{f}}(A \circ B)
 $$
  provided that 
  \begin{equation*}
\mc{A}_{\rbo} + \mc{B}_{\lbo} > 0 \text{ and } \mc{A}_{\lbo} + \mc{B}_{\rbo} > 0.
\end{equation*}
 \end{itemize}
 
 See Part I, Section 2 for more details.


\section{Resolvent kernel for asymptotically Euclidean manifolds}\label{resolvent-kernel}

In this section, we assume that $n \geq 3$, that $M$ is asymptotically Euclidean, that $P$ has no zero-resonance, and that $m' \geq 2$ if $n = 3$ (see \eqref{m'-defn}). 
We apply the same method as the case of section 3 of Part I, but with additional
difficulties due to zero modes. We recall that our operators act on
half-densities. This means that the space $L^2(M; \Omega^{1/2})=L^2(M; \Omega_b^{1/2})$ has invariant meaning (independent of a choice of metric). 

Our strategy is the same as in Part I. That is, we construct a parametrix $G(k)$ in our calculus, i.e. we want to find $G(k) \in \Psi_k^{-2; (-2,0,0), \mc{G}}(M, \Omegabht)$ that solves 
$$
(P + k^2) G(k) = \Id + E(k)
$$
with $E(k)$ `small' as $k \to 0$. As an element of the calculus,  $G(k)$ 
has a Taylor series at each boundary hypersurface of $\MMksc$. 
We use $k$ as a boundary defining function for the interior of  $\zf, \bfo, \lbo, \rbo$ and coordinates $(z,z')$ on $\zf$, $(\kappa, \kappa', y, y')$ for $\bfo$, $(z, y', \kappa')$ for $\rbo$ and $(z', y, \kappa)$ for $\lbo$. Here, on the left $M$ factor of $M^2$, $z \in M$ is written $z = (x,y)$ close to $\partial M$, where $y$ is a coordinate on $\partial M$, and $\kappa = k/x$; while primes indicate the same coordinates on the right factor. 
Using these coordinates we can write the Taylor series at the face $\textrm{f}$, for $\textrm{f} = \zf, \bfo, \rbo, \lbo$  in the form
$$
\sum_{j \geq j_0} G^j_{\textrm{f}}
$$
and this defines the coefficents $G^j_{\textrm{f}}$ uniquely. We call $G^j_{\textrm{f}}$ the `model' at order $j$ at the face $\textrm{f}$. 
Let us remark that, at the other boundary hypersurfaces of $\MMksc$, elements of the calculus have trivial expansions (i.e. are rapidly decreasing) at $\bfc, \lb, \rb$ while at $\sca$ we only need to consider the leading term (normal operator) which is well-defined independent of the choice of coordinates. 

We specify $G(k)$ by giving a finite number of models at each boundary hypersurface, together with the singularity at the diagonal, and checking compatibility of the models at adjacent faces. The parametrix can then be taken to be any element of the calculus consistent with the specified models.

\subsection{Term at $\sca$  and $\bfo$ and the singularity at the diagonal} The terms $G_{\sca}^0$, $G_{\bfo}^{-2}$ and the diagonal singularity are defined exactly as in Part I: $G_{\sca}^0$ is defined to be the inverse of $(P + k^2)_{\sca}^0$ and the diagonal singularity is defined to be the symbolic inverse of $\sigma(P + k^2)$. 
As for $G_{\bfo}^{-2}$, we recall that, in terms of coordinates $\kappa = k/x, \kappa' = k/x', y, y', k$ valid near the interior of $\bfo$, the operator $P+k^2$ reads
\begin{equation}
P+k^2=k^2\kappa^{-\ndemi-1} \Big(-(\kappa\pl_\kappa)^2+\Delta_{S^{n-1}}+(n-2)^2/4+\kappa^2+W
\Big)\kappa^{\ndemi-1}.
\label{Pbfo}\end{equation}
with respect to the flat connection on half-densities annihilating $|dg|^{1/2}$. Here $W = O(x)$ (see \eqref{pb-form}).  In terms of the  b-flat connection annihilating the b-half-density $|dg_b|^{1/2}$, our operator is to leading order at $\bfo$
\begin{equation}\label{bfo-normal-1}
\kappa^{-1}\Big(-(\kappa\pl_\kappa)^2+\Delta_{S^{n-1}}+(n-2)^2/4+\kappa^2\Big)\kappa^{-1} := \kappa^{-1} P_{\bfo} \kappa^{-1}.
\end{equation} 
The operator $P_{\bfo}$ acting on the left on $M_0\x M_0$ on b half-densities
$f(\kappa,y)|\kappa^{-1}d\kappa dy|^{\demi}$ has an inverse 
$Q_{\bfo}$ in terms of spherical harmonics $(\phi_j(y))_j$:
\begin{equation}\label{Qbfo}
\begin{split}
 Q_{\bfo}:=&\sum_{j=0}^\infty\Pi_{E_j}\Big(I_{j+\ndemi-1}(\kappa)K_{\ndemi-1+j}(\kappa')H(\kappa'-\kappa)\quad\quad\quad\quad\\
&\quad\quad\quad\quad\quad+I_{j+\ndemi-1}(\kappa')K_{\ndemi-1+j}(\kappa)H(\kappa-\kappa')\Big) \left|\frac{d\kappa dyd\kappa'dy'}{\kappa\kappa'}\right|^\demi
\end{split}\end{equation}
where here and below, $\Pi_{E_j}=\Pi_{E_j}(y,y')$ means the Schwartz kernel of the orthogonal projection on $E_j\subset L^2(S^{n-1})$, and $I_{\nu}(z),K_{\nu}(z)$ are the modified Bessel functions (see \cite{AS}).
We set
\begin{equation}
G_{\bfo}^{-2} = (\kappa \kappa') Q_{\bfo} .
\label{Gbfo}\end{equation}

The consistency of these models follows exactly as in Part I \cite{GH}, Subsection 3.5.

\subsection{Terms at zf}\label{subsect:zf}
As in Part I, the theory of b-elliptic operators given by Melrose \cite[Sec. 5.26]{APS} shows that there is a generalized inverse $Q_b$, a b-pseudo of order $-2$ for the operator $P_b$ on $L^2_b$, such that 
$$P_b Q_b = Q_b P_q = \Id - \sum_{j=0}^N\varphi_j\otimes\varphi_j.
$$
where $(\varphi_j)_j$ is a set of $L^2$ real orthonormalized (half-density) eigenfunctions of $P_b$ with eigenvalue $0$.
The normal operator $I_{\rm ff}(Q_b)$ on the front-face ${\rm ff}:=\zf\cap\bfo$ of $\zf$ is exactly as in Part I, namely the inverse of the normal operator $I_{\rm ff}(P_b) = (sD_s)^2 + (n/2-1)^{2} + \Delta_{S^{n-1}}$, given explicitly by 
\begin{equation}\label{nkernelff}
\sum_{j=0}^{\infty} \Pi_{E_j} \frac{e^{-(j+\ndemi-1) |\log s|}}{2j+n-2}|dydy'ds/s|^{\demi}
\end{equation} 
where $s:=x/x'$ is a global coordinate on the interior of ${\rm ff}$, idenfitied with $(0,\infty)_s\x S^{n-1}\x S^{n-1}$.
The $\varphi_j$ are smooth in $M$ and, by Theorem~\ref{reg}, have a polyhomogeneous expansion at the boundary $\pl M$ of the form 
\begin{equation}\label{expeigenv}
\varphi_j\sim \Big( x^{\ndemi-1} \sum_{l=0}^\infty\sum_{k=0}^la^j_{lk}(y)x^l\log^k(x) \Big) \otimes |dg_b|^{1/2}, \quad x \to 0.
\end{equation}   
The logarithmic terms in \eqref{expeigenv} complicate the parametrix construction (although not in an essential way). In order to present the construction as smoothly as possible we assume \eqref{assump1} to avoid logarithmic terms coming too early in the expansion. (See Remark~\ref{logs} for what happens if we do not make this assumption.) With this assumption, $W$ in \eqref{Pbfo} is $O(x^3)$ and we have
\begin{equation}\begin{gathered}
\varphi_j(x,y)= \sum_{i=m}^{m+2} a^j_{i0}(y)x^{\ndemi-1+i}+O(x^{\ndemi+m +2} \log x) \\ \text{ and }a^j_{i0}\in E_{i}=\ker(\Delta_{S^{n-1}}-i(n-2+i)), \quad i = 0, 1, 2. 
\end{gathered}\end{equation}
For what follows, we will write $a^j_i$ instead of $a^j_{i0}$ to avoid too many subscripts.

Since $P = x P_b x$ (see \eqref{pvspb**}), 
we obtain 
\begin{equation}\label{firstqb}
P (x^{-1}Q_b x^{-1}) =\Id -\sum_{j=0}^Nx \varphi_j \otimes x^{-1}\varphi_j.
\end{equation}
We denote by $(\psi_j)_{j=0,\dots N}$ a set of real orthonormalized eigenfunctions of $P$ for eigenvalue $0$. Then we may write
\begin{equation}
\psi_i=\sum_{j=0}^N\alpha_{ij}x^{-1}\varphi_j
\label{alpha}\end{equation}
for some matrix $(\alpha_{ij})$, whose inverse is denoted $(\alpha^{ij})$.
Decomposing in $L^2(M)$,
\[ x\varphi_j=\Pi_{\ker P}(x \varphi_j)+(1-\Pi_{\ker P})(x \varphi_j), \]
\begin{equation}
\Pi_{\ker P}( x \varphi_j):=\sum_{k=0}^N\Big(\int_{M}x \varphi_j\psi_k \Big)\psi_k
 = \sum_{l=0}^N\sum_{k=0}^N\alpha_{kl} \Big( \int_{M}\varphi_j\varphi_{l} \Big) \psi_k
  = \sum_{k=0}^N\alpha_{kj}\psi_k .
\label{magic1}\end{equation}
Thus we deduce
\begin{equation}
\sum_{j=0}^N \big( \Pi_{\ker P}(x \varphi_j) \big) \otimes x^{-1}\varphi_j=
\sum_{j=0}^N\sum_{k=0}^N\sum_{l=0}^N\alpha_{kj}\alpha^{jl}\psi_k\otimes\psi_l=
\sum_{k=0}^N\psi_k\otimes\psi_k.\label{magic2}\end{equation}
Let us denote $\psi_j^\perp:=(1-\Pi_{\ker P})(x \varphi_j)$. Then $x^{-1}\psi^\perp_j\in x^{\ndemi - 3+m-\epsilon}L_b^2(M)$  for all $\epsilon > 0$. 
Thus, for $n \geq 4$ it is in $x^{-\demi-\eps}L^2_{b}(M)$ for  $\eps>0$ (recall that $m\geq 1$ if $n=4$).
By Theorem~\ref{rit}, $P_b$ is Fredholm of index 0 on this space. 
Therefore,  $x^{-1}\psi_j^\perp$ is in ${\rm Range}(P_b)$ on $x^{-\demi-\eps}L^2_b$ if and only if it is orthogonal to the null space of $P_b$ on $x^{\demi + \epsilon} L^2(M)$. This is actually true since the null space of $P_b$ on $x^{\demi + \epsilon} L^2(M)$ is equal to the null space on $L^2(M)$ in view of the expansions \eqref{expeigenv}, and hence spanned by the $\varphi_j$. But $x^{-1} \varphi_j$ is a linear combination of the $\psi_k$ and hence $x^{-1} \varphi_j$ is orthogonal to $\psi_k^\perp$, or equivalently $\varphi_j$ is orthogonal to $x^{-1} \psi_k^\perp$. If now $n=3$ and $m\geq 2$,
then $x^{-1}\psi_j^{\perp}\in L^2_b(M)$ so we can apply the same arguments in this space instead of $x^{-\demi-\eps}L^2_b$.

This implies that there exists $\chi_j\in x^{-\demi-\eps} L^2_b(M)$ (resp. in $L^2_b(M)$) for $j=0,\dots,N$ such that
\begin{equation}\label{chi-def}
P_b\chi_j=x^{-1}\psi_j^\perp .
\end{equation}
if $n\geq 4$ (resp. if $n=3$ and $m\geq 2$).   
Assume now that $n \geq 4$ and $m' \geq 1$ for $n = 5$. (The cases $n=3,m' = 2$ and $n = 5, m' = 0$  will be discussed in Sections~\ref{Additional} and \ref{Additional5} respectively; note that $m' \geq 1$ automatically if $n = 4$.) We define our model operators at orders -2, -2, 0 and 1 at $\zf$ by
\begin{equation}\begin{aligned}
G_{\zf}^{-2} &= \sum_{j=0}^N\psi_j\otimes\psi_j \\
G_{\zf}^{-1} &= 0 \\
G_{\zf}^0 &= (xx')^{-1}\Big(Q_b+\sum_{j=0}^N(\chi_j\otimes \varphi_j+\varphi_j\otimes
\chi_j)\Big) \\
G_{\zf}^1 &= 0
\end{aligned}\label{Gzf-nullspace}\end{equation}
It is then not difficult to check that $P G_{\zf}^0 + G_{\zf}^{-2} = \Id$, which implies that $(P - k^2) G(k) - \Id$ vanishes to order $k^2$ at $\zf$ for any parametrix $G(k)$ with this expansion at $\zf$. 

We need to check consistency of these models with $G_{\bfo}^{-2}$. 
To do this we need to show that the term of order  $\rho^{-l-2}_{\zf}$ in the Taylor series for  $G_{\bfo}^{-2}$ at $\bfo \cap \zf$ agrees with the term of order $\rho_{\bfo}^{2+l}$ in the Taylor series for $G_{\zf}^l$ at $\zf \cap \bfo$, for $-2 \leq l \leq 1$, where $\rho_{\bfo} \rho_{\zf} = k$. For $l=-2$, we note that $\phi_j$ vanishes to order $n/2 - 1+m$ at $x=0$ (as a multiple of $|dg_b|^{1/2}$), hence $G_{\zf}^{-2}$ vanishes to order $ 2(n/2 - 2+m) = n - 4+2m$ at $\bfo$. Since  $n-4+2m>0$ here, we see that the restriction of this to $\zf \cap \bfo$ is zero. This matches the restriction of $G_{\bfo}^{-2}$ to $\zf \cap \bfo$ at $\zf \cap \bfo$ since $G_{\bfo}^{-2}$ vanishes to order $2$ at $\zf \cap \bfo$. Fr{o}m this we see that also $G_{\bfo}^{-2}$ and $G_{\zf}^{-1}$ are compatible. Moreover, for $\nu\in(0,\infty)\setminus\{1\}$ we have
\begin{equation}\label{inuknu}
\begin{gathered}
I_\nu(z) = \frac1{\Gamma(\nu + 1)} (\frac{z}{2})^\nu \Big( 1 +  \frac1{1 + \nu} \big( \frac{z}{2} \big)^2 + O(z^4 ) \Big)   \\
K_\nu(z) = \frac{\Gamma(\nu)}{2} (\frac{z}{2})^{-\nu} \Big( 1 +\frac{\Gamma(-\nu)}{\Gamma(\nu)}(\frac{z}{2})^{2\nu}
+ \frac1{1 - \nu} \big( \frac{z}{2} \big)^2 + O(z^4 \log z) \Big) 
\end{gathered}\qquad \text{ as } z \to 0.
\end{equation}
This implies that the next term in the Taylor series of $G_{\bfo}^{-2}$ is at order 2, so that it is also compatible with $G_{\zf}^1$.

It remains to check compatibility of $G_{\bfo}^{-2}$ and $G_{\zf}^0$. We have already observed that $G_{\bfo}^{-2}$ and $(xx')^{-1} Q_b$ are compatible. On the other hand, we know that $x^{-1} \varphi_j$ vanishes to order $n/2 - 2+m$ at $x=0$ (as a multiple of $|dg_b|^{1/2}$) and, by \eqref{chi-as}, $x^{-1}\chi$ vanishes to order $n/2 - 4+m$  at $x=0$. So the rank one terms in the third line of \eqref{Gzf-nullspace} decay to order $n-6+2m \geq -1$ at $\zf \cap \bfo$ and so restricts to zero there when multiplied by $\rho_{\bfo}^2$. This verifies compatibility of $G_{\bfo}^{-2}$ and $G_{\zf}^0$.

\subsection{Terms at $\rbo$ and $\lbo$. }
Next we construct terms $G_{\rbo}^{j}$ on $\rbo$. We will need these for three values of $j$, depending on $m$, namely $j = n/2 - 4 + i$ for $i = m', m'+1, m'+2$ where $m' = \min(m,2)$ . The terms at $\lbo$ will be determined from those at $\rbo$ by the condition that $G(k)$ has a symmetric kernel. 
We note here that in order to determine the resolvent kernel, it is not necessary to specify so many models; for this we only need the models of negative order at $\rbo$. However, the three models are required to determine the leading behaviour of the resolvent at $\rbo$ and $\lbo$ which is important for Riesz transform applications. 

We begin by determining the first few terms in the Taylor series of $G_{\zf}^{-2}$ and $G_{\zf}^0$ at $\zf \cap \rbo$. Let us begin with the kernel $Q_b$. 
Localizing near rb, the kernel of the identity vanishes identically and we have
$$
P_b Q_b = Q_b P_b =  - \sum_{j=0}^N \varphi_j \otimes \varphi_j.
$$
Using Theorem~\ref{reg} and \eqref{assump1} we can write the following asymptotic for $Q_b$ at $\rb$:
\begin{equation}\label{qbright}
Q_b= \Big( \sum_{i=0}^2\sum_{k=0}^i v_{ik}(z,y'){x'}^{\ndemi-1+i}\log^k(x') \Big) +O({x'}^{\ndemi+2}\log(x')\end{equation}
for some $v_{ik}$.
By considering the operator operating on the right variable, and using assumption \eqref{assump1}, the logarithmic terms $v_{ik}$ for $k > 0$ are absent for $i < m$, while for $i\geq m$ we have $v_{ik} = 0$ if $k \geq 2$ and $v_{i1}$ is given by 
\[v_{i1}(z,y')=\frac{1}{n-2+2i}\sum_{j=0}^Na^j_{i0}(y')\varphi_j(z)\]
since $\big((x\pl_x)^2-(n/2-1+i)^2 \big)x^{\ndemi-1+i}\log(x)=(n-2+2i)x^{\ndemi-1+i}$.
Now consider $P_b$ acting on the left variable. By matching the series 
\eqref{qbright} with the expansion \eqref{nkernelff} at $\zf \cap \bfo$ we see that for all $i$,
\begin{equation}
v_{i0}(x, y, y') = \sum_j \frac{x^{-\ndemi+1-i} \Pi_{E_i}(y,y')}{(2i+n-2)} +O(x^{-\ndemi+2-i})
\label{vi0}\end{equation}
and 
\begin{equation}
P_b v_{i0} = 0 \text{ for } i < m, \ \text{ while } \ 
P_b v_{i0} = -\sum_j \varphi_j(z) a^j_{i0}(y') \text{ for }i\geq m. 
\label{i<m}\end{equation}
Moreover, according to the assumption \eqref{assump1}, we have $v_{i0}(z,\cdot)\in E_i$ for $i\leq 2$.
Like the $a^j_i$ terms, we will denote $v_i$ instead $v_{i0}$ for simplicity of notations, and
we define $v_i:=0$ for $i<0$.

Next we consider the asymptotics of the rank one part of $G_{\zf}^0$. For this we need the asymptotics of $\chi_j$. 
The term $x^{-1}\psi_j^\perp= -\sum_{k=0}^N\alpha_{kj}x^{-1}\psi_k+\varphi_j$ has the asymptotic 
\[x^{-1}\psi_j^\perp(y)=
-\sum_{k,l=0}^N\alpha_{kj}\alpha_{kl}\sum_{i=m}^{m+2}a^l_{i}(y)x^{\ndemi-3+i}
+a_{m}^j(y)x^{\ndemi-1+m}+O(x^{\ndemi+m}\log x). \]
Theorem~\ref{reg} implies that, after setting 
\[b^j_{i}(y):=\sum_{k=0}^N\sum_{l=0}^N\alpha_{kj}\alpha_{kl}a^l_{i}(y),\]
we get the asymptotics 
\begin{equation}\begin{split}\label{chi-as}
x^{-1}\chi_j=& \sum_{i=m}^{m+2} \Big(\frac{-b^j_{i}(y)}{2(n-4+2i)}+\beta^j_{i-2}(y)\Big)x^{\ndemi-4+i} -\frac{ a_{m}^j(y)x^{\ndemi-2+m}\log x}{n-2+2m} \\ &+O(x^{n-1+m}\log^2 x).
\end{split}\end{equation}
for some $\beta^j_i\in E_{i}$ such that $\beta^j_i=0$ if $i<0$.

The leading behaviour of $G_{\zf}^{-2}$ and $G_{\zf}^0$ at $\rbo$ is therefore 
\begin{eqnarray*}
G_{\zf}^{-2}&=&x^{-1}\sum_{i=m}^{m+2}\sum_{j,k,l=0}^Na^j_{i}(y')
\alpha_{kj}\alpha_{kl}\varphi_l(z){x'}^{\ndemi-2+i}
+O({x'}^{\ndemi+1+m}\log(x'))\\
&=&x^{-1}\sum_{i=m}^{m+2}\sum_{j=0}^Nb^j_{i}(y')\varphi_j(z){x'}^{\ndemi-2+i}+O({x'}^{\ndemi+1+m}\log(x'))
\end{eqnarray*} 
and
\begin{multline}\label{Gzf0}
G_{\zf}^0=x^{-1}\sum_{j=0}^N
\left(\sum_{i=m}^{m+2}\Big(\frac{-b^j_{i}(y')}{2(n-4+2i)}+\beta_{i-2}^j(y')\Big){x'}^{\ndemi-4+i}\right) \varphi_j(z)\\
+x^{-1}\sum_{i=0}^2v_{i}(z,y'){x'}^{\ndemi-2+i}+ x^{-1}\sum_{j=0}^N  a^j_{m}(y')\chi_j(z) {x'}^{\ndemi-2+m}+O({x'}^{\ndemi-1+m'}\log^2x')\end{multline}
Note that the ${x'}^{\ndemi-2+m}\log x'$ term from $(xx')^{-1}Q_b$ cancels that of $(xx')^{-1}\sum_j\phi_j\otimes\chi_j$. 
Thus, let us define for $i = m', m'+1$
\begin{equation}
\begin{split}
G_{\rbo}^{\ndemi-4+i}:=&\frac{{\kappa'}K_{\ndemi-3+i}(\kappa')}
{\Gamma(\ndemi-3+i)2^{\ndemi-4+i}}x^{-1}\Big(\sum_{j=0}^N\beta^j_{i-2}(y')\varphi_j(z)+v_{i-2}(z,y')\Big)\\
&+\frac{{\kappa'} K_{\ndemi-1+i}(\kappa')}{\Gamma(\ndemi-1+i)2^{\ndemi+i-2}}x^{-1}\sum_{j=0}^Nb^j_{i}(y')\varphi_j(z).
\end{split}\label{grbon-4+i}
\end{equation}
Notice that $P$ annihilates these models due to \eqref{i<m}, since $l \leq i-2 \leq m-1$ in the second term. 
If $m = 0$ or $1$  we set
\begin{equation*}
\begin{split}
G_{\rbo}^{\ndemi-2+m}:=&\frac{{\kappa'}K_{\ndemi+1+m}(\kappa')}{\Gamma(\ndemi+1+m)2^{\ndemi+m}}x^{-1}\sum_{j=0}^Nb^j_{m+2}(y')\varphi_j(z)\\
 +\frac{{\kappa'}K_{\ndemi-1+m}(\kappa')}
{\Gamma(\ndemi-1+m)2^{\ndemi-2+m}}&x^{-1}\Bigg(v_{m}(z,y')
+\sum_{j=0}^N\beta^j_m(y')\varphi_j(z)+a^j_{m}(y')\chi_j(z)\Bigg).
\end{split}
\end{equation*}
Then using \eqref{i<m} and \eqref{chi-def} we compute for $m = 0,1$  
\begin{equation*}
\begin{split}
PG_{\rbo}^{\ndemi-2+m}=&\frac{{\kappa'}K_{\ndemi-1+m}(\kappa')}
{\Gamma(\ndemi-1+m)2^{\ndemi-2+m}}x \Big(P_bv_{m}(z,y')+\sum_{j}P_b\chi_j(z)a_{m}^j(y')\Big)\\
 =& \frac{{\kappa'}K_{\ndemi-1+m}(\kappa')}
{\Gamma(\ndemi-1+m)2^{\ndemi-2+m}} \Big(-\sum_{j}a_{m}^j(y')x\varphi_j(z)\\
&+\sum_{j}\big(x\varphi_j(z)-\Pi_{\ker P}(x\varphi_j(z))\big) a_{m}^j(y')\Big)\\
=&-\frac{{\kappa'}K_{\ndemi-1+m}(\kappa')}
{\Gamma(\ndemi-1+m)2^{\ndemi-2+m}}\sum_{j}b_m^{j}(y')x^{-1}\varphi_j(z)\\
=& -G_{\rbo}^{\ndemi-4+m}.
\end{split}
\end{equation*}
In this case the error term $E(k)$ will have leading behaviour at $\rbo$ no worse than $\rho_{\rbo}^{n-1+m}\log(\rho_{\rbo})$ at $\rbo$. 
When $m\geq 2$, we need to include an additional  term in $G_{\rbo}^\ndemi$ to kill the \[w(z)=w(z,y'):=v_{0}(z,y')+\sum_{j=0}^N\beta_0^j(y')\varphi_j(z)\] 
term coming from $k^2G_{\rbo}^{\ndemi-2}$ (note that this term is actually constant in $y'$ and that $\beta_0^j=0$
if $m>2$). 
 
\begin{lem}\label{complicated}
There is a half-density $\tilde v(z) \in x^{-\ndemi - 1 - \epsilon} L_b^2$ such that $P_b\til{v}(z)=x^{-2}w(z)$, with the asymptotic expansion
\begin{equation}
\til{v}(z)= \frac{x^{-\ndemi-1}}{2n(n-2){\rm Vol}(S^{n-1})}+O(x^{-\ndemi-2}), \ x \to 0. 
\end{equation}
\end{lem}

We postpone the proof of this lemma to section~\ref{proofoflemma}. Accepting this, we set
\begin{equation*}
\begin{split}
G_{\rbo}^{\ndemi}:=&\frac{{\kappa'}K_{\ndemi+3}(\kappa')}{\Gamma(\ndemi+3)2^{\ndemi+2}}x^{-1}\sum_{j=0}^Nb^j_{4}(y')\varphi_j(z)\\
& +\frac{{\kappa'}K_{\ndemi+1}(\kappa')}
{\Gamma(\ndemi+1)2^{\ndemi}}x^{-1}\Bigg(v_{2}(z,y')
+\sum_{j=0}^N\beta_2^j(y')\varphi_j(z)+a^j_{2}(y')\chi_j(z)\Bigg) \\
&+ \frac{{\kappa'} K_{\ndemi-1}(\kappa')}
{\Gamma(\ndemi-1)2^{\ndemi-2}}{x'}^{n-2}x^{-1}\til{v}(z).
\end{split}
\end{equation*}
and then again $PG_{\rbo}^{\ndemi} = -G_{\rbo}^{\ndemi - 2}$, so the error term at $\rbo$ has leading term at order $\rho_{\rbo}^{\ndemi+1}\log(\rho_{\rbo})$.

We next observe that the dependence of $G_{\rbo}^{\ndemi - 4 + i}(z, \kappa', y')$ on $\kappa'$ and $y'$, for $i = m', m'+1$ is always of the form $\kappa' K_{n/2 - 1 + j}(\kappa') b_j(y')$ where $b_j \in E_j$. It follows that $P$ acting in the right variable kills these models, or equivalently that $P G_{\lbo}^{\ndemi - 4 + i} = 0$ for $i = m', m'+1,m'+2$.
In addition, we automatically gain two orders at the left boundary since $P = x P_b x$ and $x=0$ at $\lbo$.  Therefore the error term has leading order $n/2 +1 + m'$ at the left boundary. 

Using the formulae in (\ref{inuknu}),
one easily checks that $G_{\rbo}^{\ndemi-4+i}$ matches with $G_{\zf}^{j}$ for $i=m',m'+1,m'+2$ and $j\in\{-2,-1,0,1\}$ at $\zf\cap\rbo$. The  $G_{\rbo}^{\ndemi-4+i}$ also match with $G_{\bfo}^{-2}$, using the asymptotics of $v_{l0}(z)$ and of $I_\nu(\kappa)$.  Also note that when $m' = 2$ the matching of $G_{\rbo}^{\ndemi}$ and $G_{\bfo}^{-2}$ involves the subleading, rather than leading, term of $I_{\ndemi - 1}(\kappa)$ at $\kappa = 0$ in $G_{\bfo}^{-2}$.

\subsection{Additional term when $n=3$ and $m\geq 2$}\label{Additional}
In this case, all terms we have constructed above are kept the same except
the term $G_{\zf}^1$. Indeed, to 
match $G_{\zf}^1$ with $G_{\bfo}^{-2}$, we see that the term of order $z^{\demi}$ in the asymptotic
\[K_{\demi}(z)=\frac{\Gamma(\demi)}{2}(\frac{z}{2})^{-\demi}(1-z+O(z^2))\]
implies the asymptotic of $G_{\bfo}^{-2}$ at $\zf$
\[G_{\bfo}^{-2}=\kappa\kappa' (I_{\ff}(Q_b)-k(xx')^{-\demi}\phi_0(y)\phi_0(y')+O(\rho_{\zf}^2))\]
and we thus set the term
\[G_{\zf}^{1}=-(xx')^{-1}{\rm Vol}(S^2)v_{0}(z)v_{0}(z')\]
which is compatible with $G_{\bfo}^{-2}$ in view of (\ref{vi0}) and $\phi_0(y)=({\rm Vol}(S^2))^{-\demi}$. 
As for matching of $G_{\zf}^1$ with $G_{\rbo}^{-\demi}$, it comes from the second asymptotic term 
in the $\kappa'K_{1/2}(\kappa')v_{0}(z)$ term of $G_{\rbo}^{-\demi}$.  

\subsection{The case $n=5$ with $m=0$}\label{Additional5}

In this case, we see that all our parametrix can be constructed similarly, except that a term at order 
$1$ at $\zf$ needs to be added to match with $G_{\rbo}^{-3/2}$:  indeed $z^\nu K_{\nu}(z)$ has a cubic term in its expansion only for $n=5$, which explains why we need a nonzero $G_{\zf}^1$ term, and this term will force us to add a $G_{\zf}^{-1}$
as well. 

For this part we assume that $M$ has one end to simplify, then $\dim E_0=1$  and we can always suppose that $\psi_0$ is a normalized eigenfunction decaying like $c_0x^{1/2}|dg_b|^{\demi}$ 
with $c_0:=\sum_{j=0}^N\alpha_{0j}a^j_0$ and that the other $\psi_j$ are in $x^{1/2+\eps}L^2_b(M)$ for small $\eps>0$ and orthogonal to $\psi_0$. Then using Theorem~\ref{rit} we can see that there exists $\theta$ such that
$P\theta=\psi_0$ with $\theta\sim e_0x^{-5/2}|dg_b|^\demi$ for some constant $e_0$. Indeed, $x^{-1}\psi_0$
is in the range of $P_b$ acting in $x^{-3/2-\eps}L^2_b$ since it is orthogonal to the Null space of $P_b$ on $x^{3/2+\eps}L^2_b$ (this null space being spanned by the $(x\psi_j)_{j>0}$), thus there exist $\til{\theta}\in x^{-3/2-\eps}L^2_b$ with $P_b\til{\theta}=x^{-1}\psi_0$ and
$\theta:=x^{-1}\til{\theta}$ satisfies $P\theta=\psi_0$.
Fr{o}m Theorem \ref{reg} and the equation $P_bx^{-1/2}=2x^{-\demi}+O(x^{\demi})$, the asymptotic of $\til{\theta}$ is given by  
\[\til{\theta}= e_0 x^{-3/2}+ \demi c_0x^{-1/2}+O(x^{1/2}), \quad e_0\in\cc\]
and by considering Green's formula on $\lim_{\eps\to 0}\int_{x>\eps}(P_b\til{\theta})x\psi_0-\til{\theta}P_b(x\psi_0)=1$ we find $e_0=-(3 c_0{\rm Vol}(S^4))^{-1}$.
Then we define 
\begin{equation}\label{casen=5m=0}
G_{\zf}^{1}:=-c_0^2 {\rm Vol}(S^4)\Big(\theta \otimes \psi_0 +\psi_0\otimes\theta\Big), \quad G_{\zf}^{-1}:=
c_0^2{\rm Vol}(S^4)\psi_0\otimes\psi_0.\end{equation}
so that $PG_{\zf}^1=-G_{\zf}^{-1}$.
The term of order $\rho_{\zf}^1$ in $k^{-3/2}G_{\rbo}^{-3/2}$ comes from the $\kappa'^{3/2}$ coefficient in the expansion
of $K_{3/2}(\kappa')$, that is $k{x'}^{-5/2}/3$ times
\[\sum_{j=0}^Nb_0^jx^{-1}\varphi_j(z)=\sum_{j,k,l=0}^N\alpha_{kj}x^{-1}\varphi_j(z)\alpha_{kl}a_{0}^l=\sum_{j,l=0}^N\alpha_{0j}x^{-1}\varphi_j(z)\alpha_{0l}a_0^l=c_0\psi_0(z),
\] 
which implies consistency between $G_{\zf}^1, G_{\zf}^{-1}$ and $G_{\rbo}^{-3/2}$. Note that $G_{\zf}^1, G_{\zf}^{-1}$  also match with $G_{\rbo}^{-1/2}$ if we modify
$G_{\rbo}^{-1/2}$ by adding the term 
\[\frac{\kappa'K_{3/2}(\kappa')}{\Gamma(3/2)2^{1/2}}\frac{c_0^3}{{\rm Vol}(S^4)}\psi_0(z),
\] 
as is straightforward to check using (\ref{inuknu}). 

\begin{remark} The analysis in this Subsection is only necessary if we wants to specify the leading term at $\rbo$ and the $k^{-1}$ term at $\zf$. If one is content to construct a parametrix with an error that iterates away, then one can specify $G_{\zf}^j$ for $j = -2, -1, 0$ as in \eqref{Gzf-nullspace} and leave $G_{\zf}^1$ unspecified. 
\end{remark}

\begin{remark}\label{log6} A similar phenomenon occurs in dimension 4 and 6: the log term in $K_2(z)$ forces a term at order $(2,1)$, i.e. at order $k^2 \log k$ at $\zf$, which requires a nonzero $G_{\zf}^{0,1}$, i.e. at order $\log k$. We shall not emphasize this here although it is important in that it contributes to the leading order behaviour of the spectral measure at $k = 0$. 
\end{remark}

\subsection{Error term and resolvent} 
The error term $E$ defined by $(P + k^2) G = \Id + E$ now vanishes to order $2$ at $\zf$, order $1$ at $\bfo$ and $\sca$, and order $n/2 - 1 + m'$ at $\rbo$ and $n/2 + 1 + m'$ at $\lbo$ (possibly with log terms at leading order at $\lbo$ and $\rbo$, and at $\zf$ for $n=4,6$ only --- see Remark~\ref{log6}).  Therefore it iterates away in the sense of Remark 2.12 of Part I, and as above the inverse $\Id + S = (\Id + E)^{-1}$ exists for small $k$ and lies in the calculus. Here $S$ has index sets contained in $(2,1) \cup (2 + \mc{N}')$ at  $\zf$, $1 + \mc{N}'$ at $\bfo$, $1$ at $\sca$, $(n/2+1+m',1) \cup n/2 + 1 +m' + \mc{N}'$ at $\lbo$, $(n/2 - 1 + m',1) \cup n/2 - 1 +m' + \mc{N}'$ at $\rbo$ for some nonnegative integral index set $\mc{N}'$ and empty at $\bfc$, $\lb$ and $\rb$. 
The resolvent itself is given by $R(k) = G(k) + G(k) S(k)$ and satisfies

\begin{theo}\label{thm:nullspace} Assume that $n \geq 5$, or that $n = 3$ or $4$ with the additional condition that $P$ has no zero-resonances. We also assume that $P$ satisfies $\eqref{assump1}$. 
Let $m'$ be defined as in \eqref{m'-defn} and assume \eqref{m'-cond}. 
Then for small $k_0$, $k < k_0$, the resolvent $R(k) = (P + k^2)^{-1}$ on half-densities satisfies
\begin{equation}
R(k) \in \Psi^{-2, (-2, 0, 0), \mc{R}}(M, \Omegab)
\end{equation}
where the index family $\mc{R}$ satisfies 
\begin{align*}
\mc{R}_{\zf} &\subset (-2,0) \cup (-1,0) \cup (0,1) \cup \mc{N}'', \\ 
\mc{R}_{\bfo} &\subset -2 + \mc{N}'', \\
\mc{R}_{\sca} &= 0  \\
\mc{R}_{\lbo} = \mc{R}_{\rbo} &\subset n/2 - 4 + m' + \mc{N}'' 
\end{align*} for some 
integral (with log-terms) index set $\mc{N}'' \geq 0$, and the empty set at all other faces. Moreover, the leading terms of the resolvent $R(k)$ at $\zf, \bfo, \sca, \lbo, \rbo$ are equal to the leading terms of the parametrix $G(k)$ as defined above. We also have $R_{\zf}^{-1} = G_{\zf}^{-1} = 0$ in all cases except when $n=5,m=0$ where $R_{\zf}^{-1}=G_{\zf}^{-1}$ is given by $\eqref{casen=5m=0}$; $R_{\zf}^{0,1} = 0$ for $n \neq 4,6$; and $R_{\zf}^0$ is always equal to $G_{\zf}^0$ modulo a finite rank term with values in $\ker_{L^2} P$. 
\end{theo}

\begin{remark} In the case that $M$ is flat $\RR^n$,  Jensen \cite{Jensen} obtains $R_{\zf}^{-2} = P_0$, the projection onto the zero eigenspace, $R_{\zf}^{-1} = 0$ for $n \geq 6$ and gives the expression $P_0 V G_3 V P_0$ for $R_{\zf}^{-1}$ when $n=5$, where $G_3$ is the operator with kernel constant and equal to $1/3$. 
Note that 
$$
\int_{\RR^5} V \psi_j = - \int_{\RR^5} \Delta \psi_j = \lim_{R \to \infty} \int_{\partial B(R,0)} \partial_r \psi_j(R, \omega) d\omega 
$$
is equal to $c_0 {\rm Vol}(S^4)$ if $j=0$ and $0$ otherwise. Using this and ${\rm Vol}(S^4) = 8\pi^2/3$ we can check that this is the same as our $R_{\zf}^{-1}$ for $n=5$. The agreement on $R_{\zf}^j$ for $j = -2$ and $j = -1, n\geq 6$ is clear,  so our results on the singular part of the resolvent at $\zf$ are in agreement with Jensen's.  
\end{remark}

\begin{remark}\label{logs} If we do not make assumption \eqref{assump1}, our proof  allows to construct a weaker parametrix (in particular with only the first term at $\rbo,\lbo$) which shows, with the composition formula \cite[Prop 2.10]{GH}, that the same result holds if $n>6$ but with another index set $\mc{R}$ which has same properties than $\mc{R}$ except that $\mc{R}_{\zf}\subset (-2,0)\cup \mc{R'}$, and we can not specify the leading terms of $R(k)$ at $\rbo,\lbo$.
In general the proof could be well adapted, there will be additional logarithmic terms at $\rbo$ and $\lbo$ that need to be specified, which would make the paper extremely technical. If $m=0$, then there may be terms of order 
$\rho_{\rbo}^{\ndemi - 3} \log \rho_{\rbo}$, $\rho_{\rbo}^{\ndemi - 2} (\log \rho_{\rbo})^2$ and $\rho_{\rbo}^{\ndemi - 2} \log \rho_{\rbo}$. The extra models will all be similar in structure to $G_{\ndemi - 4 + i}$ for $i = 1, 2$ and create no essential difficulties. 
\end{remark} 

\subsection{Proof of Lemma~\ref{complicated}}\label{proofoflemma}
To prove the existence of $\tilde v$, we use Theorem~\ref{rit}, which tells us that $x^{-2}w(z)$ is  in the range of $P_b$ acting on $x^{-\ndemi - 1 - \epsilon} L^2_b(M)$ if and only if $x^{-2}w(z)$ is orthogonal to the null space of $P_b$ on $x^{\ndemi + 1 + \epsilon} L_b^2(M)$. 

Let us  
begin with the equation $P G_{\zf}^0 = \Id - \sum_j \psi_j \otimes \psi_j$, which implies that $G^0_{\zf} \psi_j$ is in the span of the $\psi_k$. If $m \geq 2$, the $\psi_k$ vanish to order at least $n/2$ at $x=0$ (as a b half-density). Now consider a $\psi\in\ker P$ that lies in $x^{\ndemi + \epsilon} L^2$ (and hence 
actually decreasing to order $\ndemi + 1$), corresponding to
a $\varphi_j \in \ker P_b$ that lies in $x^{\ndemi + 1 +
\epsilon} L^2$. Consider applying $G_{\zf}^0$ to such a
$\psi$. Fr{o}m the assumption $m \geq 2$ and the fact that this
lies in the span of the $\psi_j$ we see that the result is
$O(x^{\ndemi - 2})$ as a b-half-density. However, the kernel
of $G_{\zf}^0$ vanishes to order $-2$ at $\zf$ and $n/2 - 2$
at $\lbo$, while the lift of $\psi$ to $\zf$ (in the right
factor) vanishes to order $n/2 + 1$ at $\rbo$ and $\lbo$. 
Fr{o}m this, we see that $G_{\zf}^0 \psi$ apparently only vanishes to order $n/2 - 2$ at $x=0$, with a coefficient  obtained from the 
integral of $\psi$ against the value of $x^{-\ndemi + 2} G_{\zf}^0$ restricted to $\zf \cap \lbo$. But $x^{-\ndemi + 2} G_{\zf}^0$ restricted to $\zf \cap \lbo$ is, from \eqref{Gzf0} (with the left and right variables switched) equal to 
$$
w(z')-\sum_{j=0}^N\frac{b_2^j(y)}{2n}\varphi_j(z').
$$
We conclude that this integrated against $\psi(z')$ vanishes:
\begin{equation}\label{ifm=2}
x^{\ndemi-2}\int_{M}{x'}^{-1}\Big(w(z')-\sum_{j=0}^N\frac{b_2^j(y)}{2n}\varphi_j(z')\Big)\psi(z')=0 
\end{equation}
We also note that the coefficient of $x^{\ndemi - 2}$ must be constant in $y$. However $b_2^j \in E_2$ which is orthogonal to constants so integrating in $y$ (and recalling $w$ is independent of $y$ we deduce that 
$$
\int_{M}{x'}^{-1} w(z') \psi(z') = 0 
$$
as claimed.

So choose $\tilde v \in x^{-\ndemi - 1 - \epsilon} L^2_b(M)$ with $P_b \tilde v = x^{-2} w$, which is
defined modulo elements in the null space of $P_b$ 
in $x^{-n/2-1-\eps}L^2_b(M)$. Then Theorem~\ref{reg} and the asymptotic behaviour of $w(z)$ using (\ref{vi0}) 
show that
\[
\til{v}(z)= \frac{x^{-\ndemi-1}}{2n(n-2){\rm Vol}(S^{n-1})}+x^{-\ndemi-1}\beta(y)+O(x^{-\ndemi-2})\]
for some $\beta\in E_2$ (recall that an $L^2$-normalized element of $E_0$ is $({\rm Vol}(S^{n-1}))^{-\demi}$).
It is necessary, in order to match with $G_{\bfo}^{-2}$, that $\beta=0$, and we will see that this can be obtained 
by adding a term in the null space of $P_b$ in $x^{-\ndemi-1-\eps}L^2_b(M)$. By Proposition~\ref{prop:complementary}, this is possible if and only if $\beta$ is orthogonal to the subspace $G_2 \subset E_2$ consisting of the leading asymptotics of elements of $\ker P_b$ that are $\sim x^{\ndemi+1}$ as $x \to 0$, or equivalently of elements $\psi$ of $\ker P$ that are $\sim x^{\ndemi}$ as $x \to 0$.

We now prove this. Let $\psi$ be an element of $\ker P$ that is equal to $x^{\ndemi} \gamma(y) + O(x^{\ndemi + 1})$, as a b-half-density, 
with $\gamma \in E_2$. Note by Green's formula we have
\begin{multline*} \int_{S^{n-1}}\beta(y)\gamma(y)dy = \int_{S^{n-1}}\gamma(y)(\beta(y)+c_0)= \frac1{n-2} 
\lim_{\eps\to 0}\int_{x>\eps}P_b\til{v}(z)x\psi(z) \\
= \frac1{n-2} 
\lim_{\eps\to 0}\int_{x>\eps} x^{-1} w(z)\psi(z)
\end{multline*}
so it is enough to show that 
\begin{equation}
\lim_{\eps\to 0}\int_{x>\eps} x^{-1} w(z) \psi(z) = 0.
\label{want}\end{equation}
Again we look at $G_{\zf}^0 \psi$. Since now we are assuming $\psi \sim x^{\ndemi}$, Melrose's Pushfoward theorem \cite{cocdmc} tells us that 
$$
G_{\zf}^0 \psi \sim x^{\ndemi - 2} \log x d_0(y) + x^{\ndemi - 2} d_1(y) + O(x^{\ndemi - 1} \log x),
$$
where $d_1$ is given by a sum of terms 
\[ \lim_{\eps\to 0}\int_{x'>\eps}{x'}^{-1}\Big(w(z')-\sum_{j=0}^N\frac{b_2^j(y)}{2n}\varphi_j(z')\Big)\psi(z')\]
\[+ \int_{s \leq \epsilon^{-1}} I_{\rm ff}(Q_b)(s,y,y')(s)^{\ndemi-1}\gamma(y') ds dy' \]
and $d_0 = d_1 = 0$ since $m \geq 2$. 
To compute the second term, we use the formula \eqref{nkernelff} to obtain the limit $\gamma(y)/2n$. For the first term, we can take 
$\psi=x^{-1}\varphi_j$ without loss of generality, whence $\gamma = a^j_{20}(y)$, and we compute 
\[\sum_{l=0}^Nb_2^l(y)\int_M{x'}^{-1}\varphi_l(z'){x'}^{-1}\varphi_j(z')=a^j_{2}(y)\]
since $\int_M {x}^{-1}\varphi_l(z){x}^{-1}\varphi_j=\sum_{k=0}^N\alpha^{kl}\alpha^{kj}$
and $b_2^l(y)=\sum_{k=0}^N\sum_{i=0}^N\alpha_{il}\alpha_{ik}a_{2}^k(y)$. This proves that \eqref{want}, and completes the proof of the lemma. 

\section{Resolvent kernel for asymptotically conic manifolds}
We make the same assumptions as in the previous section, but allow $M$ to be asymptotically conic rather than Euclidean. It is quite similar
to the asymptotically Euclidean case but more terms come in
the parametrix. We recall the definitions of $m$ and $m'$ in \eqref{m'-defn}
 (now $m$ is not necessarily an integer), define $\nu_m = n/2-1+m \in N_\partial$, and assume \eqref{m'-cond}. Then the eigenfunctions have asymptotics 
\[\varphi_j=\sum_{\nu_m\leq \nu_i\leq \nu_m+2} a^j_{\nu_i}(y)x^{\nu_i}+O(x^{\ndemi+1+\eps})\]
for some $\eps>0$. 
Now the term $Q_b$ has asymptotic on $\rb$
\[Q_b=\sum_{k=0}^1\sum_{\nu_0\leq\nu_j\leq \nu_0+2}v_{\nu_j,k}(z,y'){x'}^{\nu_j}(\log x')^k+ O({x'}^{\nu_0+2+\eps})\]
for some $\eps>0$, with $v_{\nu_j,k}=0$ if $k=1,j<m$, and 
\[P_bv_{\nu_i,0}(z,y')=\left\{\begin{array}{ll}
0 & \text{ if }i<m\\
-\sum_{j=0}^N a_{\nu_i}(y')\varphi_j(z) & \text{ if }i\geq m
\end{array}\right.
\] 
with $v_{\nu_j}(z,\cdot):=v_{\nu_j,0}(z,\cdot)\in E_{\nu_j}$ and $v_{\nu_j}(z,y')$ 
having leading behaviour (obtained from $I_{\ff}(Q_b)$) 
\[v_{\nu_j}(x,y,y')= \frac{\Pi_{E_j}(y,y')}{2\nu_j}x^{-\nu_j}+O(x^{-\nu_{j-1}}\log x).\]
 We also define $v_z:=0$ if  $z\notin N_{\pl}$ to match 
with our following notations.

The $\chi_j$ are defined as in the previous section, and  have asymptotics 
\begin{equation*}
x^{-1}\chi_j= \!\!\!\!\!\!\! \sum_{\nu_m\leq \nu_i\leq \nu_m+2} \!\!\!\! \Big(\frac{b^{j}_{\nu_i}(y)}{4(-\nu_i+1)}+\beta^j_{\nu_i-2}(y)\Big)x^{\nu_i-3} -\frac{ a_{\nu_m}^j(y)x^{\nu_m-1}\log x}{2\nu_m} +O(x^{\nu_m-1+\eps}).
\end{equation*}
for some $\beta^j_{z}\in E_{z}$ such that $\beta^j_z=0$ for $z<\nu_0$. 

We set as before $G_{\bfo}^{-2}:=\kappa\kappa'Q_{\bfo}$ with 
\begin{multline}\label{bfo-conic-model}
 Q_{\bfo}:= \sum_{j=0}^\infty\Pi_{E_j}\Big(I_{\nu_j}(\kappa)K_{\nu_j}(\kappa')H(\kappa'-\kappa)+
I_{\nu_j}(\kappa')K_{\nu_j}(\kappa)H(\kappa-\kappa')\Big) \\
\left|\frac{d\kappa dyd\kappa'dy'}{\kappa\kappa'}\right|^\demi
\end{multline}
and (except for $n=5,m'=0$ -- see Section~\ref{n5m0-conic}) define
\begin{equation}\begin{aligned}
G_{\zf}^{-2} &= \sum_{j=0}^N\psi_j\otimes\psi_j \\
G_{\zf}^{-1} &= 0 \\
G_{\zf}^0 &= (xx')^{-1}\Big(Q_b+\sum_{j=0}^N(\chi_j\otimes \varphi_j+\varphi_j\otimes
\chi_j)\Big) \\
G_{\zf}^1 &= 0 
\end{aligned}\end{equation}

The leading behaviour of $G_{\zf}^{-2}$ and $G_{\zf}^0$ at $\rbo$ is therefore 
\begin{eqnarray*}
G_{\zf}^{-2}&=&x^{-1}\sum_{\nu_m\leq \nu_i\leq \nu_m+2}\sum_{j,k,l=0}^Na^j_{\nu_i}(y')
\alpha_{kj}\alpha_{kl}\varphi_l(z){x'}^{\nu_i-1}
+O({x'}^{\nu_m+1+\eps})\\
&=&x^{-1}\sum_{\nu_m\leq\nu_i\leq \nu_m+2}\sum_{j=0}^Nb^j_{\nu_i}(y')\varphi_j(z){x'}^{\nu_i-1}+O({x'}^{\nu_m+1+\eps})
\end{eqnarray*} 
and
\begin{multline}
G_{\zf}^0=x^{-1}\sum_{j=0}^N
\left(\sum_{\nu_m\leq\nu_i\leq\nu_m+2}
\Big(\frac{-b^j_{\nu_i}(y')}{4(\nu_i-1)}+\beta_{\nu_i-2}^j(y')\Big){x'}^{\nu_i-3}\right) \varphi_j(z)\\
+x^{-1}\sum_{\nu_0\leq\nu_i\leq \nu_0+2}v_{\nu_i}(z,y'){x'}^{\nu_i-1}+ x^{-1}\sum_{j=0}^N  a^j_{\nu_m}(y')\chi_j(z) {x'}^{\nu_m-1}+O({x'}^{\ndemi-1+m'+\eps}).\end{multline}
Note that the ${x'}^{\nu_m-1}\log x'$ term from $(xx')^{-1}Q_b$ cancels that of $(xx')^{-1}\sum_j\phi_j\otimes\chi_j$. 
Thus for $i$ such that $\nu_0+m'\leq \nu_i< \nu_0+2+m'$, we define
\begin{equation}\label{leadingrbo}
\begin{split}
G_{\rbo}^{\nu_i-3}:=&
\frac{{\kappa'}K_{\nu_i-2}(\kappa')}
{\Gamma(\nu_i-2)2^{\nu_i-3}}x^{-1}\Big(\sum_{j=0}^N\beta^j_{\nu_i-2}(y')\varphi_j(z)+v_{\nu_i-2}(z,y')\Big)\\
&+\frac{{\kappa'} K_{\nu_i}(\kappa')}{\Gamma(\nu_i)2^{\nu_i-1}}x^{-1}\sum_{j=0}^Nb^j_{\nu_i}(y')\varphi_j(z).
\end{split}
\end{equation}
If $m'<2$, we set
\begin{equation*}
\begin{split}
G_{\rbo}^{\nu_m-1}:=&\frac{{\kappa'}K_{\nu_m}(\kappa')}
{\Gamma(\nu_m)2^{\nu_m-1}}x^{-1}\Bigg(v_{\nu_m}(z,y')
+\sum_{j=0}^N\beta^j_{\nu_m}(y')\varphi_j(z)+a^j_{\nu_m}(y')\chi_j(z)\Bigg) \\
 &+\frac{{\kappa'}K_{\nu_m+2}(\kappa')}{\Gamma(\nu_m+2)2^{\nu_m+1}}x^{-1}\sum_{j=0}^Nb^j_{\nu_m+2}(y')\varphi_j(z)
\end{split}
\end{equation*}
so that $PG_{\rbo}^{\nu_m-1}=G_{\rbo}^{\nu_m-3}$. If now $m'=2$, we set 
\begin{equation*}
\begin{split}
G_{\rbo}^{\nu_0+1}:=&\frac{{\kappa'}K_{\ndemi+3}(\kappa')}{\Gamma(\ndemi+3)2^{\ndemi+2}}x^{-1}\sum_{j=0}^Nb^j_{\nu_0+4}(y')\varphi_j(z)\\
& +\frac{{\kappa'}K_{\ndemi+1}(\kappa')}
{\Gamma(\ndemi+1)2^{\ndemi}}x^{-1}\Bigg(v_{\nu_0+2}(z,y')
+\sum_{j=0}^N\beta_{\nu_0+2}^j(y')\varphi_j(z)+a^j_{\nu_0+2}(y')\chi_j(z)\Bigg) \\
&+ \frac{{\kappa'} K_{\ndemi-1}(\kappa')}
{\Gamma(\ndemi-1)2^{\ndemi-2}}{x'}^{n-2}x^{-1}\til{v}(z,y').
\end{split}
\end{equation*}
where $\til{v}$ is defined like in Lemma \ref{complicated} (thus actually $\til{v}(z,y')$ is constant in $y'$) so that 
\[P_b\til{v}(z,y')=v_{\nu_0}(z,y')+\sum_{j=0}^N\beta_{\nu_0}^j(y')\varphi_j(z).\]
Then again $PG_{\rbo}^{\ndemi} = -G_{\rbo}^{\ndemi - 2}$, so the error term at $\rbo$ has leading term at order $\rho_{\rbo}^{\nu_0+m'-1+\eps})$ for some $\eps>0$. 

\subsection{Additional term when $n=3$}
Now when $n=3$ and $m>1$, there is an additional term at $\zf$ at order $1$ as when $\pl M=S^{2}$, this is
\[G_{\zf}^1:=-(xx')^{-1}{\rm Vol}_{h_0}(\pl M)v_0(z)v_0(z')\]
and it is there to match with the second term in the asymptotic of $G_{\bfo}^{-2}$ at $\zf$. 

\subsection{Terms at $\zf$ when $n=5$, $m' = 0$}\label{n5m0-conic}
This works just as in the Euclidean case; we define $G_{\zf}^{\pm1}$ by \eqref{casen=5m=0}, replacing ${\rm Vol}(S^4)$ with ${\rm Vol}(\partial M)$ (with respect to the metric $h(0)$ from \eqref{metricconic}). 

\subsection{Error term and resolvent}
The error term $E$ defined by $(P + k^2) G = \Id + E$ vanishes to order $1+\eps$ at $\zf$ (for some $\eps>0$ depending on a finite number of $\nu_j$), at order $1$ at $\bfo$ and $\sca$, and order $n/2 - 2 + m' + \epsilon$ at $\rbo$ and $n/2  + m' + \epsilon$ at $\lbo$.  Then it iterates away in the sense of Remark 2.12 of Part I, and the inverse $\Id + S = (\Id + E)^{-1}$ exists for small $k>0$ and lies in the calculus, where $S$ has index sets included 
in $1+ \mc{Z}$ at  $\zf$, $(1,0) \cup (1+\mc{Z})$ at $\bfo$, $1$ at $\sca$, $n/2+m'+ \mc{Z}$ at $\lbo$, $n/2-2+m'+ \mc{Z}$ at $\rbo$ and empty at $\bfc$, $\lb$ and $\rb$ for some index set $\mc{Z}$ such that $\mc{Z}\geq \eps$ for some $\eps>0$. 
The resolvent itself is given by $R(k) = G(k) + G(k) S(k)$ and satisfies

\begin{theo}\label{res-conic-1} Assume that $n \geq 5$, or that $n = 3$ or $4$ with the additional condition that $P$ 
has no zero-resonances. We assume that $P$ satisfies $\eqref{assump1}$ and $\eqref{m'-cond}$. 
Then for small $k_0$, $k < k_0$, the resolvent $R(k) = (P + k^2)^{-1}$ on half-densities satisfies
\begin{equation}
R(k) \in \Psi^{-2, (-2, 0, 0), \mc{R}}(M, \Omegab)
\end{equation}
where the index family $\mc{R}$ satisfies 
\begin{align*}
\mc{R}_{\zf} &\subset (-2,0) \cup (-1,0) \cup (-1+\mc{Z}'), \\ 
\mc{R}_{\bfo} &\subset (-2,0) \cup (-1,0) \cup(-1 + \mc{Z}'), \\
\mc{R}_{\sca} &= 0  \\
\mc{R}_{\lbo} = \mc{R}_{\rbo} &\subset n/2 - 4 + m' + \mc{Z}' 
\end{align*} for some index set $\mc{
Z}' \geq \eps$ where $\eps>0$ depends on $N_\pl$, and the index set at the other faces are empty. Moreover, the leading terms of the resolvent $R(k)$ at $\zf, \bfo, \sca, \lbo, \rbo$ are equal to the leading terms of the parametrix $G(k)$ as defined above. We also have $R_{\zf}^{-1} = G_{\zf}^{-1}$ (which vanishes except in the case $n=5,m'=0$).  
\end{theo}


\section{Resolvent: asymptotically Euclidean manifolds, dimension 3}\label{dim3}
In this section we suppose that $n=3$ and that $P$ has $L^2$ kernel and possibly also a zero-resonance.
We use similar notations and same method as in Section~\ref{resolvent-kernel} since this is almost exactly the same problem, the only difference being that $m'$ may be $0$ or $1$. To simplify (to avoid log-terms coming too early in expansions) 
we make assumption \eqref{assump1} throughout this section.

\subsection{Terms $G_{\zf}^{-2}$, $G_{\zf}^0$, $G_{\rbo}^0$}
Transposed into the b-problem as before we have an 
orthonormal basis $(\varphi_j)_{j=0,\dots,N}$ of $\ker P_b$. 
To simplify exposition, we assume that $M$ has only one end. If $P$ has a resonance at $0$, then
one of the $\varphi_j\in\ker P_b$ is $\sim ax^{\demi} |dg_b|^{\demi}$ at $x=0$
for some $a\in\cc$, which we may suppose is $\varphi_0$. We sometimes denote it
$\varphi$ instead to point out the difference with the other eigenstates. We
shall use the convention  $\varphi_0=0$ if $0$ is not a resonance.  
There also may be zero-modes $\varphi_j \sim x^{3/2} a^j_1$ at $x=0$, with $a^j_1 \in \ker (\Delta_{S^2} - 2)$ (these were ruled out in Section~\ref{resolvent-kernel} by our assumption that $m \geq 2$ when $n = 3$.) 

We can decompose in an orthogonal sum
\[\ker P_b=(\ker P_b\cap x^{3/2+\eps}L^\infty(M))\oplus \mc{H}\oplus \cc\varphi\]
for some $\mc{H}$ of dimension $d\leq 3=\dim\ker (\Delta_{S^2}-2)$. 
Then we can suppose that $(\varphi_i)_{i=1,\dots,d}$ is an orthonormal 
basis of $\mc{H}$ and in general
\[\varphi_j= \Big( a_{1}^j(y)x^{3/2}+a^j_{2}(y)x^{5/2}+O(x^{\frac{7}{2}}\log x) \Big) |dg_b|^{1/2} \]
\[{\rm with }\quad a_{i}^j\in \ker (\Delta_{S^2}-2),\quad a^{j}_{2}\in\ker(\Delta_{S^2}-6)\]
for any $1 \leq j\leq N$ with $a_1^j=0$ if $j>d$. The mapping 
$\varphi_j \to a_1^j$
identifies $\mc{H}$ with a subspace (still noted $\mc{H}$) of $\ker(\Delta_{S^2}-2)$ and we denote by $\mc{H}^\perp$
an orthogonal complement of $\mc{H}$ in $\ker(\Delta_{S^2}-2)$. Let  now
 $(\phi_1,\phi_2,\phi_3)$ an orthonormal basis of $\ker (\Delta_{S^2}-2)$, and we
note $a_{1}^j=\sum_{l=1}^da^{j}_{1l}\phi_l$ for some $a^j_{1l}\in\cc$. 

There exists $Q_b$ such that $P_bQ_b=Q_bP_b={\rm Id}-\sum_{j=1}^N\varphi_j\otimes\varphi_j$
as in Section~\ref{resolvent-kernel}. Localizing near $\rb$, the identity becomes $P_bQ_b=Q_bP_b=-\sum_j\varphi_j\otimes\varphi_j$ and considering $P_b$ acting on the right factor $z'$, this gives the following asymptotic for $Q_b$ at $\rb$
\[Q_b(z,z')=\left\{\begin{array}{ll}
v_{0}(z,y'){x'}^{\demi}+O({x'}^{3/2}) & \textrm{ if }\varphi=0 \\
-a\varphi(z){x'}^\demi\log x'+v_1(z,y'){x'}^{\demi}+O({x'}^{3/2})& \textrm{ if }\varphi\not=0
\end{array}\right.\]
for some $v_{i}(z,y')=v_{i}(z)$, $i=0,1,2,$ constant in $y'$, which, by considering the operator acting on the left variable $z$, must satisfy 
\begin{equation}\label{pbvo}
P_bv_{0}(z)=0, \quad P_bv_1(z)=-a\varphi(z)
\end{equation} 
This existence of $v_0$ is guaranteed by Proposition~\ref{prop:complementary} if $P$ has no zero-resonance. 
Similarly, if $\varphi\not=0,$ the existence of an element $u\in x^{-1/2}L^2_b$ such that $P_bu=\varphi$
is a consequence of Theorem \ref{rit}: indeed, $P_b$ is Fredholm in $x^{-1/2-\eps}L^2_b$ with index $1$ and
$\varphi$ is in the range of $P_b$ in this space since it is orthogonal to the kernel of $P_b$ in $x^{1/2+\eps}L^2_b$.
Notice that $u$ is uniquely determined modulo $\ker P_b$, it is polyhomogeneous, and by adding a constant times $\varphi$ to $u$, it can be chosen to have asymptotic 
\[u=bx^{-\demi}+ax^\demi\log x +O(x^{3/2}\log x)\]
for some $b\in\cc$. Then $v_1(z)=-au(z)+\beta\varphi(z)+\eta(z)$ for some $\beta\in\cc$ and $\eta\in\ker P_b\cap x^{3/2}L^\infty(M)$.
The normal operator at the front face $I_{\rm ff}(Q_b)$ given in
(\ref{nkernelff}) implies that $v_0(z)$ has the aymptotic 
\[v_{0}(z)=\frac{1}{{\rm vol}(S^{2})} x^{-\demi}+O(x^{\demi})\]
when $\varphi=0$. While when $\varphi\not=0$, Green formula with the asymptotics of $u,\varphi$ give the relation
\[\cjg P_bu,\varphi\cjd=1=\lim_{\eps\to 0}\int_{x=\eps}((x\pl_xu) \varphi -u (x\pl_x\varphi)=-ab{\rm  Vol}(S^2)=-4\pi ab\]
thus $ab=-(4\pi)^{-1}$. 
We now define 
\[\til{Q}_b:=Q_b+u\otimes\varphi+\varphi\otimes u\]
which satisfy $P_b\til{Q}_b=1-\sum_{j=1}^N\varphi_j\otimes\varphi_j$.
We obtain
\[P \big( x^{-1}\til{Q}_bx^{-1} \big)=\Id-\sum_{j=1}^Nx \varphi_j\otimes x^{-1}\varphi_j.\]
Gathering this information, we find that, when $\varphi\not\equiv 0$, the asymptotic behaviour of $\til{Q}_b$ at $\rb$ is given by 
\begin{equation}\label{asymptilQ}
\til{Q}_b(z,z')=b\varphi(z){x'}^{-\demi}+(\beta\varphi(z)+\eta(z)){x'}^\demi+O({x'}^{3/2}\log x').
\end{equation}
We remark that the ${x'}^\demi\log x'$ term of $Q_b$ at $\rb$ is cancelled by the $\varphi\otimes u$ 
term of $\til{Q}_b$.
 
As before, we denote by  $(\psi_j)_{j=1,\dots N}$ 
an orthonormal basis of eigenfunctions of $P$ and express these in terms of the $\varphi_j$ by \eqref{alpha}. Following the argument  in Section~\ref{resolvent-kernel} (see \eqref{magic1}), we have 
\[\sum_{j=1}^N\Pi_{\ker P}(x\varphi_j)\otimes x^{-1}\varphi_j=
\sum_{j=1}^N\sum_{k=1}^N\sum_{l=1}^N\alpha_{kj}\alpha^{jl}\psi_k\otimes\psi_l=
\sum_{k=1}^N\psi_k\otimes\psi_k.\]
We proceed exactly as before: we denote $\psi_j^\perp:=(1-\Pi_{\ker P})(x\varphi_j)$, then $x^{-1}\psi^\perp_j=O(x^{-\demi})$
is not anymore in $L_b^2$ but instead is in $x^{-\demi-\eps}L_b^2$. The extension of $P_b$
to $x^{-\demi-\eps}L_b^2$ is Fredholm by Theorem~\ref{rit}.  
Hence $x^{-1}\psi_j^\perp$ is in ${\rm Range}(P_b)$ on $x^{-\demi-\eps}L_b^2$
if $x^{-1}\psi_j^\perp \perp \ker_{x^{\demi+\eps}L_b^2}P_b=\ker_{L_b^2}P_b$. This is satisfied 
since
\[0=\int_{M}\psi^\perp_j \bar{\psi}_k =\int_M x^{-1}\psi_j^\perp x\bar{\psi}_k,\] 
and $(x\psi_j)_j$ is a basis of $\ker_{x^{\demi+\eps}L^2} P_b$.
This implies that there exists $\chi_j\in x^{-\demi-\eps} L^2_b$ for $j=0,\dots,N$ such that
\[P_b\chi_j=x^{-1}\psi_j^\perp.\]
As in Section~\ref{resolvent-kernel} we set the term at order $-2$ at $\zf$ to 
\[
G_{\zf}^{-2}:=\sum_{j=1}^N\psi_j\otimes\psi_j,\quad
G_{\zf}^0:=(xx')^{-1}\Big(\til{Q}_b+\sum_{j=1}^N(\chi_j\otimes \varphi_j+\varphi_j\otimes
\chi_j)\Big)
\]
The term $x^{-1}\psi_j^\perp= -\sum_{k=1}^N{\alpha_{kj}}x^{-1}\psi_k+\varphi_j$ has the asymptotic 
\[x^{-1}\psi_j^\perp(y)=
-\sum_{k,l=1}^N{\alpha_{kj}}\alpha_{kl}(a_1^l(y)x^{-\demi}+a_{2}^l(y)x^{\demi})+O(x^{\frac{3}{2}}\log x)\]
which implies that, after setting 
\begin{equation}
b_i^j(y):=\sum_{m=1}^N\sum_{l=1}^N{\alpha_{mj}}\alpha_{ml}a_i^l(y),
\label{bijy}\end{equation}
we get the asymptotic
\begin{equation}
x^{-1}\chi_j(x,y)=x^{-3/2}(-\demi b_1^j(y)+\beta^j_{-1})+x^{-1/2}(-\frac{1}{6}b_{2}^j(y)+\beta^j_{0})+O(x^{\demi}\log x)
\label{chij-asympt}\end{equation}
for some $\beta^j_{0},\beta^j_{-1}\in \cc$.
It is important to notice that $\chi_j$ is determined modulo elements in the null space of $P_b$
in $x^{-1/2-\eps}L^2_b$, but this space contains either the function $v_0\sim (4\pi)^{-1}x^{-1/2}$ (if $0$ is not resonance)
or a function $\varphi\sim ax^{\demi}$ (if $0$ is resonance), which mean in first case that any $\beta^j_{-1}$ can be supposed to be $0$ by adding a constant times $v_0$ to $\chi_j$ while in the second case any $\beta^j_0$ can be taken to be $0$ by adding a constant times $\varphi$. 
  
Observe that $G_{\zf}^{-2}$ is $O(\rho_{\bfo})$ and $x^{-1}\chi_j\otimes x^{-1}\varphi_j,x^{-1}\varphi_j\otimes x^{-1}\chi_j$ are $O(\rho_{\bfo}^{-1})$. 

Note that we will often use the identity
\begin{equation}\label{identity}
\sum_{j=1}^Nx^{-1}\varphi_j(z){b^j_{i}(y')}=\sum_{l=1}^N\Pi_{\ker P}(x\varphi_l)(z){a^l_{i}(y')}
\end{equation}
which follows from \eqref{alpha}, \eqref{magic1} and \eqref{bijy}. 

\subsection{Term $G_{\bfo}^{-2}$}
We next define $G_{\bfo}^{-2}$. As above and Part I, we use coordinates $(\kappa=\frac{k}{x}, \kappa'=\frac{k}{x'}, y, y', k)$ which are valid near the interior of $\bfo$. We write $G_{\bfo}^{-2} = (\kappa \kappa') Q_{\bfo}$ where $Q_{\bfo}$ solves 
$$
P_{\bfo} Q_{\bfo} = \Id
$$
and $P_{\bfo}$ is given by \eqref{bfo-normal-1}. 
We can write the general solution $Q_{\bfo}$ in terms of spherical harmonics but we'll need
an additional finite rank term when $\varphi\not\equiv 0$, in order to match with the $u\otimes\varphi+\varphi\otimes u$
term of $G_{\zf}^0$: we set ($H=$Heaviside)
\[ Q_{\bfo} =\sum_{j=0}^\infty \Pi_{E_j} \Big(I_{j+\demi}(\kappa)K_{\demi+j}(\kappa')H(\kappa'-\kappa)+
I_{j+\demi}(\kappa')K_{\demi+j}(\kappa)H(\kappa-\kappa')\Big)\]\[
+ c_{0} \Pi_{E_0} K_{\demi}(\kappa)K_{\demi}(\kappa'), \quad\textrm{ with } E_{0}(y,y')=1/\textrm{Vol}(S^2)=\frac1{4\pi} \]
for some $c_{0}\in\cc$ to determine. In higher dimensions, this additional term is too singular at $\kappa = \kappa' = 0$, i.e. at $\bfo \cap \zf$ and therefore cannot appear. This coefficient $c_0$ is set to be $0$ when $\varphi\equiv 0$ while 
when $\varphi\not\equiv 0$, it is constructed to match with $G_{\zf}^0$: we find, using asymptotics for $K_{\demi}(z)$ given by
\begin{equation}\label{k1/2}
K_{\demi}(z) = \sqrt{\frac{\pi}{2}} \Big( z^{-\demi} - z^{\demi} + O(z^{\frac{3}{2}}) \Big), \quad z \to 0,
\end{equation}
that $c_{0} = 2/\pi$. 
Indeed, using leading behaviours of $u$ and $\varphi$ when $\varphi\not=0$, the asymptotic of $G_\zf^0$ near the corner $\zf\cap\bfo$ is  
\[G_{\zf}^0= (xx')^{-1} \Big( I_{\rm ff}(Q_b)+ab \big( x^\demi{x'}^{-\demi} + {x'}^\demi{x}^{-\demi} \big) +O(\rho^{-1}_{\bfo}) \Big) \]
whereas the asymptotic
of $G_{\bfo}^{-2}$ near this corner is
\[G_{\bfo}^{-2} = k^2(x x')^{-1} \bigg(  I_{\rm ff}(Q_b)+ \frac{c_0}{8} \Big( k^{-1}(xx')^{\demi} -  ({x'}^{\demi}x^{-\demi} + {x}^{\demi}{x'}^{-\demi})
+O(\rho_{\zf}) \bigg).\]
Recalling that $ab=-1/4\pi$ from our analysis of $G_\zf^0$, we find that we must set $c_0 = 2/\pi$ (we emphasize that it is the \emph{subleading} term and not the leading term of the $K_{\demi}(\kappa) K_{\demi}(\kappa')$ term in $G_{\bfo}^{-2}$ that matches with $u \otimes \phi + \phi \otimes u$). The \emph{leading} term will be at order $-1$ at $\zf$  and it will require us to have a corresponding nonzero $G_{\zf}^{-1}$.  Note that, with such a choice for $c_0$, the term of order $k$ in the expansion of $G_{\bfo}^{-2}$ at $\zf$ is given by
\begin{equation}\label{thirdterm} 
(8\pi)^{-1}k(x^{-5/2}{x'}^{-\demi}+{x'}^{-5/2}x^{-\demi}),
\end{equation}
which will be useful to define $G_{\zf}^1$.

\subsection{Term $G_{\zf}^{-1}$ when $0$ is resonance, first attempt}
To match with the term at order $-1$ of $G_{\bfo}^{-2}$ at $\zf$, we need a term of the form  $\psi\otimes\psi$
where $\psi$ is a zero-resonance. This has the form $\psi = cx^{-1}\varphi$ for some $c\not=0$ to determine, plus an element of the $L^2$ null space. But matching with $G_{\bfo}^{-2}$ implies that
the leading term of $\psi$ is 
\begin{equation}
\psi(z)=\sqrt{\frac1{4 \pi x}} +O(x^{1/2})
\label{psi-norm}\end{equation}
where we chose the sign to be $+$ by convention.
This means $ca=(4 \pi)^{-\demi}$; the leading behaviour of $\psi$ then agrees with the canonical resonance of Jensen-Kato \cite{JK}. 
As for the $L^2$ null space part (which does not affect the leading behaviour), it is natural to ask that $\psi$ is orthogonal, in a generalized sense, to the $\psi_k$. Observe that for $1 \leq j \leq d$,  $\psi \cdot \psi_j$ has asymptotic $a^j_1(y) x^3 + O(x^4)$ as $x \to 0$ and is therefore not $L^1$; however, since $a_j^1$ is orthogonal to constants, the integral of $a^j_1(y)$ in $y$ vanishes and therefore 
$$
\lim_{\epsilon \to 0} \int_{M \cap \{ x \geq \epsilon \}} \psi \cdot \psi_k \ \text{ exists.}
$$
We denote this limit $\langle \psi, \psi_k \rangle$. We choose $\psi:=cx^{-1}\varphi-\sum_{k=1}^N\cjg cx^{-1}\varphi,\psi_k\cjd\psi_k$ to be the unique zero-resonance satisfying \eqref{psi-norm} and with $\langle \psi, \psi_k \rangle = 0$ for $1 \leq j \leq N$. (If there is no zero-resonance then we take $\psi = 0$ in the formulae below.) This agrees with the canonical zero-resonance of Jensen-Kato in the case that $(M,g) = (\RR^3, \delta)$; see Section~\ref{ss:JK} for more details. 
We now temporarily define  (the tilde indicates that this will be modified below)
\[
\tilde G_{\zf}^{-1}:=\psi \otimes \psi
\]
which matches with $G_{\bfo}^{-2}$. 

\subsection{Term $G_{\rbo}^{-3/2}$}
We construct $G_{\rbo}^{-3/2}$ like in Section~\ref{resolvent-kernel} with $m=1$ if $\varphi\equiv 0$, while if $\varphi\not\equiv 0$, we shall add a term. In any case, we define it to be
\begin{equation}\label{grbo-3/2}
\begin{gathered}
G_{\rbo}^{-3/2}:=
\frac{{\kappa'} K_{3/2}(\kappa')}{\Gamma(3/2)2^{1/2}}x^{-1}\sum_{j=1}^Nb^j_{i}(y')\varphi_j(z)
+ \frac{\kappa'K_{1/2}(\kappa')}{\Gamma(1/2)2^{-\demi}}(4\pi)^{-\demi}\psi(z).
\end{gathered}\end{equation}
By construction and from the asymptotics of $K_{1/2},K_{3/2}$, this terms matches with $G_{\zf}^{-2},G
_{\zf}^{-1},G_{\bfo}^{-2}$ and also with $G_{\zf}^0$ when $\varphi\equiv 0$. The consistency between this term and $G_{\zf}^0$ when $\varphi\not\equiv 0$ is less straightforward but it is satisfied using (i) from Lemma~\ref{lemtechnic},  the expansion
\begin{equation}\label{asymptkk1/2}
\frac{\kappa'K_{1/2}(\kappa')}{\Gamma(1/2)2^{-\demi}}=\kappa'^{\demi}-\kappa'^{3/2}+\demi{\kappa'}^{5/2}+O({\kappa'}^{7/2})
\end{equation}  
and $-(4\pi)^{-1/2}c=b$ in view of the  identities $ab=-(4\pi)^{-1}$ and $ca=(4\pi)^{-1/2}$.
 
\subsection{The terms $G_{\zf}^{1}, G_{\zf}^{-1}$}

The term at order $1$ at $\zf$ has to match with $G_{\bfo}^{-2}$ and $G_{\rbo}^{-3/2}$. 
We first assume there is a resonance, i.e. $\varphi\not\equiv 0$.
Then $x^{-1}\psi$ is orthogonal to any $x\psi_k$, thus to the null space of $P_b$ in $x^{3/2+\eps}L^2_b$
for $\eps\in(0,1)$, so we deduce that it is in the range of $P_b$ acting on $x^{-3/2-\eps}L^2_b$.  Therefore there exists
$\chi\in x^{-3/2}L^\infty(M)$ such that
\[P_b \chi=x^{-1}\psi \textrm{ or equivalently }P(x^{-1}\chi)=\psi.\]
Moreover, considering the leading asymptotic of $\psi$, we must have  
\begin{equation}\label{gamma}
\chi =-\demi(4\pi)^{-\demi}x^{-3/2}+\gamma(y)x^{-3/2} +O(x^{-\demi})\end{equation}
for some $\gamma\in E_1$ since $(-(x\pl_x)^2+1/4)x^{-3/2}=-2x^{-3/2}$. Let us assume that $\gamma=0$, which can
be arranged, possibly after adding a term $x^{-1}v$ for some $v\sim -\gamma x^{-3/2}$ in the null space of $P_b$ 
in $x^{-3/2-\eps}L^2_b$. So as not to interrupt the exposition we defer the proof to  Lemma \ref{lemtechnic}.  
We then see from (\ref{thirdterm}) that a term 
\begin{equation}\label{gzf1firsttry}
-\Big((x^{-1}\chi)\otimes \psi+\psi \otimes (x^{-1}\chi)\Big) \end{equation}
of order $1$ at $\zf$ would match with $G_{\bfo}^{-2}$, while from (\ref{asymptkk1/2}) and (\ref{grbo-3/2}) it is clear that
it would match with the $K_{1/2}(\kappa')$ terms of $G_{\rbo}^{-3/2}$ (i.e. the first line of (\ref{grbo-3/2})) at $\rbo\cap \zf$.

It remains to add a part at $\zf$, order $1$, that will match the $K_{3/2}(\kappa')$ terms of $G_{\rbo}^{-3/2}$.
These terms have leading asymptotic at zf, modulo $O(\rho_{\zf}^{7/2})$
\begin{equation*}
\begin{split}
G_{\rbo}^{-3/2}=& \big({\kappa'}^{-1/2}-\demi {\kappa'}^{3/2} +\frac{1}{3}{\kappa'}^{5/2} \big)
x^{-1}\sum_{p=1}^3\sum_{l=1}^d\sum_{m,j=1}^N\alpha_{mj}{\alpha_{ml}a_{1p}^l\phi_p(y')}\varphi_j(z)\\
=& \big({\kappa'}^{-1/2}-\demi {\kappa'}^{3/2} +\frac{1}{3}{\kappa'}^{5/2} \big)
\sum_{p=1}^3\sum_{l=1}^d\Pi_{\ker P}(x\varphi_l)(z){a_{1p}^l\phi_p(y')}\\
\end{split}
\end{equation*}
In particular, we see that a term at order $1$ at $\zf$ must be added to (\ref{gzf1firsttry}).
Fr{o}m Theorem~\ref{rit} and Proposition~\ref{prop:complementary},  there exists $\theta_j\in x^{-\frac{3}{2}}L^\infty(M)$ for $j=1,2,3$ such that 
\[P_b\theta_j=\left\{\begin{array}{ll}
\varphi_j & \textrm{ for }  j\leq d\\
0 & \textrm{ for }j>d
\end{array}\right. \text{ with } \] 
\[\theta_j=\sum_{l=1}^3c^j_1(y)x^{-\frac{3}{2}}+O(x^{-\demi}), \quad c^j_1(y)=\sum_{l=1}^3c^j_{1l}\phi_l(y)\in\ker(\Delta_{S^2}-2)\]
for some $c^j_{1l}\in\cc$. Then $(c^j_1)_{j=d+1,\dots,3}$ form a basis of $\mc{H}^\perp$.
Fr{o}m Green's formula applied to 
$$
\delta_{ij} = \int_{M \cap \{ x \geq \epsilon \} } \Big( \langle P_b \theta_i, \varphi_j \rangle - \langle \theta_i, P_b \varphi_j \rangle \Big), \quad j \leq d, \ i \leq 3,  
$$
we obtain, by letting $\eps\to 0$,
 \[ 3\int_{S^2}c^i_1(y)a^j_1(y)dy=  \delta_{ij},  \quad j\leq d, i\leq 3.\]
Let $C$ be the $3\x 3$ matrix whose entries are $C_{jl}:=c^j_{1l}$, 
$A$ the $3\x 3$ matrix whose entries are $A_{jl}:=a^j_{1l}$ for $j\leq d$ and $A_{jl}:=c^j_{1l}/3$ for $j>d$ (thus invertible), 
we have $C= \frac1{3} (A^T)^{-1}$.
Note that $\theta_j$ is determined up to a $O(x^{-1/2})$ term and since this coefficient of order is constant, 
then, if $\varphi\equiv 0$, we can add a constant times $v_0$ so that $\theta_j(x,y)=\sum_{l=1}^3c^j_{1l}\phi_l(y)+O(x^\demi)$.
In the case $\varphi\not\equiv 0$, the term of order $x^{-1/2}$  vanishes automatically;  see Lemma \ref{lemtechnic} below.

Let us define $\til{G}^1_{\zf}$, (the tilde indicates that this will be modified below):
\begin{multline*} \tilde G_{\zf}^{1}:=\Big(\sum_{i=1}^d\sum_{j=1}^3d_{ij}x^{-1} {\theta_j\otimes\Pi_{\ker P}(x\varphi_i)}  +
d_{ij}\Pi_{\ker P}(x \varphi_i)\otimes{x'}^{-1}\theta_j\Big) \\ +(xx')^{-1}ev_0\otimes v_0, 
\end{multline*}
when $\varphi\equiv 0$, while if $\varphi\not\equiv 0$
\begin{multline*}
\tilde G_{\zf}^{1}:=\Big(\sum_{i=1}^d\sum_{j=1}^3d_{ij}x^{-1} {\theta_j\otimes\Pi_{\ker P}(x\varphi_i)}  +
d_{ij}\Pi_{\ker P}(x \varphi_i)\otimes{x'}^{-1}\theta_j\Big) \\
-\Big((x^{-1}\chi)\otimes \psi+\psi\otimes (x^{-1}\chi)\Big)
\end{multline*}
where $d_{ij},e\in\cc$  are parameters to determine and $D:=(d_{ij})$ is the correponding $ d\x 3$ matrix. To match with $G_{\rbo}^{-3/2}$ and using the fact
that $(\phi_l)_{l=1,\dots,d}$ and $(\Pi_{\ker P}(x\varphi_l))_{l=1,\dots,N}$ are respectively linearly independent,  
this forces to have as matrices $DC= \frac1{3} \Pi_dA$ where $\Pi_d:M_{3\x 3}(\cc)\to M_{d\x 3}(\cc)$ is the canonical projection from $3\x 3$ matrices to $d\x 3$ matrices (so that $\Pi_d A$ is the matrix $A$ with the $c^j_{1l}$ changed to $0$). We conclude that 
\begin{equation}
D= \Pi_dAA^T.
\label{D}\end{equation} 
Note that $(d_{ij})_{1\leq i,j\leq d}=\cjg a_1^i,a^j_1\cjd_{L^2(S^2)}$ is a symmetric $d\x d$ matrix that depends only on 
$(\varphi_j)_{1\leq j\leq d}$.

Consider the matching conditions $\til{G}_{\zf}^1 \leftrightarrow G_{\bfo}^{-2}$. The terms in $\tilde G_{\zf}^1$ involving $\theta_j$ are $O(\rho_{\bfo}^{-2})$ so have no effect at order $-2$ on $\bfo$. When $\varphi\equiv0$, the remaining term involving $v_0$ is of order $\rho_{\bfo}^{-3}$ and so has to match with the $\rho_{\zf}^{3}$ coefficient of $G_{\bfo}^{-2}$. 
Using \eqref{inuknu} we see that this comes only from the $j=0$ term of 
\eqref{Qbfo} (since only for $\nu = 1/2$ does the second term in the expansion of $K_\nu$ differ from the leading term by one order). The terms match provided we choose $e=-{\rm vol}(S^2)=-4\pi$. If, on the contrary $\varphi\not\equiv 0$, then 
we have seen that the terms involving $\chi,\psi$ are those who match with $G_{\bfo}^{-2}$, 
so this proves consistency between $\til{G}_{\zf}^{1}$ and $G_{\bfo}^{-2}$.

Now we have on $\zf$ 
\[P \tilde{G}^1_{\zf}=-\psi\otimes\psi+
\Big(\sum_{j=1}^d\sum_{i=1}^d d_{ij} (x\varphi_i)\otimes\Pi_{\ker P}(x \varphi_j)\Big).\] 
This is not quite what we want since we require $PG_{\zf}^{1}=-G_{\zf}^{-1}$ with  $PG_{\zf}^{-1}=0$ in order to construct
a parametrix at order $2$ at $\zf$. To remedy this we  decompose $x\varphi_j=\Pi_{\ker P}(x\varphi_j)+\psi_j^\perp$ and then modify $\tilde{G}_{\zf}^1$ to
\begin{equation*}
\begin{split}
G_{\zf}^1:=&\sum_{i,j=1}^d \Big( d_{ij}x^{-1}{(\theta_j-\chi_j)\otimes\Pi_{\ker P}(x \varphi_i)}+
d_{ij}{x'}^{-1}\Pi_{\ker P}(x\varphi_i)\otimes(\theta_j-\chi_j)\Big)\\
&+ \Big(\sum_{i=1}^d\sum_{j=d+1}^3d_{ij}x^{-1}\theta_j\otimes\Pi_{\ker P}(x\varphi_i)+
d_{ij}{x'}^{-1}\Pi_{\ker P}(x\varphi_i)\otimes\theta_j\Big)\\
&-4 \pi (xx')^{-1} v_0\otimes v_0
\end{split}
\end{equation*}
when $P$ has no resonance at $0$ while 
\begin{equation*}
\begin{split}
G_{\zf}^1:=&\sum_{i,j=1}^d \Big( d_{ij}x^{-1}{(\theta_j-\chi_j)\otimes\Pi_{\ker P}(x \varphi_i)}+
d_{ij}{x'}^{-1}\Pi_{\ker P}(x\varphi_i)\otimes(\theta_j-\chi_j)\Big)\\
&+ \Big(\sum_{i=1}^d\sum_{j=d+1}^3d_{ij}x^{-1}\theta_j\otimes\Pi_{\ker P}(x\varphi_i)+
d_{ij}{x'}^{-1}\Pi_{\ker P}(x\varphi_i)\otimes\theta_j\Big)\\
&-\Big((x^{-1}\chi)\otimes \psi+\psi\otimes (x^{-1}\chi)\Big)
\end{split}
\end{equation*}
when $0$ is resonance.
 so that  
\[PG_{\zf}^1=\sum_{i,j=1}^dd_{ij}{\Pi_{\ker P}(x\varphi_i)\otimes\Pi_{\ker P}(x\varphi_j)} -\psi\otimes\psi\]
Then we modify $G_{\zf}^{-1}$ to set 
\begin{equation}
G_{\zf}^{-1}:=-\sum_{i,j=1}^dd_{ij}{\Pi_{\ker P}(x\varphi_i)\otimes\Pi_{\ker P}(x\varphi_j)}+\psi\otimes\psi
\label{Gzf-1}\end{equation}
which is a symmetric expression satisfying $PG_{\zf}^{-1}=0$ and matching with $G_{\zf}^1$ in the sense
$PG_{\zf}^1=-G_{\zf}^{-1}$. (We also observe that $G_{\zf}^{-1}$ agrees with that in Section~\ref{Additional} when $m=2$ since then $d=0$, hence $D = 0$, and $\psi=\chi=0$.)

Let us now prove several deferred results in the following lemma:

\begin{lem}\label{lemtechnic}

(i) The coefficients $\beta^j_{-1}$ from \eqref{chij-asympt} satisfy
\begin{equation}\label{simplify}
\sum_{j=1}^N\beta_{-1}^jx^{-1}\varphi_j=(4\pi)^{-\demi}\sum_{k=1}^N\cjg cx^{-1}\varphi,\psi_k\cjd \psi_k. \end{equation}
ii) The term $\gamma$ in the asymptotic (\ref{gamma}) of $\chi$ is zero if $\dim\mc{H}=3$ or can made $0$ by possibly adding an element in $\mc{H}$.\\
iii) If $\varphi\not\equiv 0$, the coefficient of $x^{-1/2}$ in the expansion of $\theta_j$ is $0$.
\end{lem}
\textsl{Proof}: the proof of all claims follow from Green's formula. To prove (i), we apply Green's formula to 
$$
\lim_{\epsilon \to 0 } \int_{M \cap \{ x \geq \epsilon \} } \langle P_b\chi_j,c\varphi \rangle - \langle \chi_j,P_b(c\varphi) \rangle  
$$
and obtain 
$
(4\pi)^{\demi} \beta^j_{-1}  = \langle \Pi_{\ker P} (x \varphi_j), cx^{-1}\varphi \rangle ;
$
using \eqref{magic1} and the identity $x^{-1}\varphi_j=\sum_{i=1}^N\alpha^{ji}\psi_i$ we obtain \eqref{simplify}. 

For ii), we apply Green formula  
\[\int_{x>\eps} (P_b\chi)x\psi_k -\chi P_b(x\psi_k)=\int_{x=\eps}(x\pl_x\chi).x\psi_k-\chi.x\pl_x(x\psi_k)\]
for any $\psi_k\in\ker P$ and let $\eps\to 0$. The limit has to be zero by construction of $\chi$ and since $P_b(x\psi_k)=0$,
but from the expansions of $\chi$, (\ref{expeigenv}) and (\ref{alpha}), we see that the right hand side has for limit (up to non zero constant) $\sum_{j=1}^d\alpha_{kj}\int_{S^2}\gamma(y)a_{1}^j(y)$. Now $\alpha_{kj}$ is invertible so 
$\gamma$ is orthogonal to the functions $(a_1^j)_{j=1,\dots,d}$ on $S^2$. If $d=3$, we get $\gamma=0$, otherwise from Proposition \ref{prop:complementary}, there is a function with the asymptotic $\sim x^{-3/2}\gamma(y)$ in the null space  
of $P_b$ thus i) is proved. As for iii), this is similar by using Green formula on $\int_{x>\eps}
 (P_b\theta_j)\varphi -(P_b\varphi)\theta_j$ and taking the limit as $\eps\to 0$ (which gives $0$).
\qed

\subsection{The term $G_{\rbo}^{-1/2}$}
To see what $G_{\rbo}^{-1/2}$ should be, we consider the asymptotics of $G_{\zf}^j$ at $\rbo$. These are first (where we use \eqref{identity})
\begin{equation*}\begin{aligned}
G_{\zf}^{-2}=& {x'}^{1/2} \sum_{j=1}^d\Pi_{\ker P}(x\varphi_j)(z)a^j_1(y')+{x'}^{3/2}\sum
_{l=0}^N\Pi_{\ker P}(x\varphi_l)(z){a_{2}^l(y')} \\
&+O({x'}^{5/2} \log x') \\
G_{\zf}^{-1} =& -{x'}^{1/2} \sum_{i,j=1}^d d_{ij}\Pi_{\ker P}(x\varphi_i)(z)b^j_1(y')-(4\pi)^{-1/2}{x'}^{-1/2}\psi(z)\\
&+{x'}^{1/2}\phi(y')\psi(z)+O({x'}^{3/2}) \\
\end{aligned}
\end{equation*}
where $\phi\in E_1$ is the term in the expansion $\psi(z)=(4\pi x)^{-1/2}+\phi(y)x^{1/2}+O(x^{3/2})$, 
and if there is no resonance    
\begin{equation*}\begin{aligned}
G_{\zf}^0 &= -\demi {x'}^{-3/2} \sum_{j=1}^d\Pi_{\ker P}(x\varphi_j)(z){a_1^j(y')}\\ 
&+ {x'}^{-1/2} \Big(x^{-1}v_{0}(z)+\sum_{j=0}^N\Pi_{\ker P}(x\varphi_j)(-\frac{1}{6}a_{2}^j(y')+\beta_0^j)\Big)+O({x'}^{1/2})\\
G_{\zf}^1 &= \frac{1}{3} {x'}^{-5/2}  \sum_{l=1}^d\Pi_{\ker P}(x\varphi_l)(z)a_{1}^l(y')\\
& - {x'}^{-3/2} \Big( x^{-1}v_0(z)- \demi\sum
_{i,j=1}^dd_{ij}\Pi_{\ker P}(x\varphi_i) b^j_1(y') \Big) +O({x'}^{-1/2}).
\end{aligned}\end{equation*}
while if there is a resonance
\begin{equation*}\begin{aligned}
G_{\zf}^0=& -\demi {x'}^{-3/2} \sum_{j=1}^d\Pi_{\ker P}(x\varphi_j)(z){a_1^j(y')}\\ 
& + {x'}^{-1/2} \Big(\beta x^{-1}\varphi(z)+x^{-1}\eta(z)-\frac{1}{6}\sum_{j=0}^N\Pi_{\ker P}(x\varphi_j){a_{2}^j(y')}\Big)+O({x'}^{1/2})\\
 G_{\zf}^1=&\frac{1}{3} {x'}^{-5/2}  \sum_{j=1}^d\Pi_{\ker P}(x\varphi_j)(z)a_{1}^j(y')\\
&- {x'}^{-3/2} \Big(\mu\psi(z)- \sum
_{i,j=1}^dd_{ij}\Pi_{\ker P}(x\varphi_i)(\demi b^j_1(y')-\beta_{-1}^j) \Big) +O({x'}^{-1/2}). 
\end{aligned}\end{equation*}
where $\mu\in\cc$ is defined from the expansion $\chi=-\demi(4\pi)^{-1/2}x^{-3/2}+\mu x^{-1/2}+O(x^{1/2})$ and $\eta,\beta$ are defined in the construction of $G_{\zf}^0$ from the expansion (\ref{asymptilQ}) of $\til{Q}_b$ at $\rbo$.
We need to check that the coefficient of ${x'}^{-1/2 - j}$ in the expansion of $G_{\zf}^j$ agrees with the coefficient of ${\kappa'}^{j+1/2}$ in the expansion of $G_{\rbo}^{-1/2}$, and finally $\beta_{-1}^j,\beta_0^j$ come from the expansion of $\chi_j$. Let $W\in C^\infty([0,1])$ a function with support in $[0,1/2)$ which is equal to $1$ near $0$. We set 
\begin{equation*}
\begin{split}
G_{\rbo}^{-1/2}:=&
x^{-1}v_{0}(z)\frac{{\kappa'} K_{\demi}(\kappa')}{\Gamma(1/2)2^{-1/2}}+W(\kappa')\sum_{j=1}^N\beta_0^j\Pi_{\ker P}(x\varphi_j)\\
&-\frac{{\kappa'} K_{\frac{3}{2}}(\kappa')}{\Gamma(3/2)2^{1/2}}\sum_{i,j=1}^3 d_{ij} \Pi_{\ker P}(x\varphi_i)(z)b^j_1(y')\\
&+\frac{{\kappa'}K_{\frac{5}{2}}(\kappa')}{\Gamma(5/2)2^{3/2}}
\sum_{l=0}^N\Pi_{\ker P}(x\varphi_l){a_{2}^l(y')}\end{split}
\end{equation*}
when $\varphi\equiv 0$ or 
\begin{equation*}
\begin{split}
G_{\rbo}^{-1/2}:=&-{\kappa'}^{3/2}W(\kappa')\Big(\mu\psi(z)+\sum_{i,j=1}^d\beta_{-1}^jd_{ij}\Pi_{\ker P}(x\varphi_j)\Big)\\
&+x^{-1}(\beta\varphi(z)+\eta(z))\frac{{\kappa'} K_{\demi}(\kappa')}{\Gamma(1/2)2^{-1/2}}\\
&-\frac{{\kappa'} K_{\frac{3}{2}}(\kappa')}{\Gamma(3/2)2^{1/2}}\sum_{i,j=1}^3 d_{ij} \Pi_{\ker P}(x\varphi_i)(z)b^j_1(y')\\
&+\frac{{\kappa'}K_{\frac{5}{2}}(\kappa')}{\Gamma(5/2)2^{3/2}}
\sum_{l=0}^N\Pi_{\ker P}(x\varphi_l){a_{2}^l(y')}\end{split}
\end{equation*}if $\varphi\not\equiv 0$.
By construction, both satisfy these matching conditions with $G_{\zf}^i$, $i=-2,\dots,1$. 
The matching $G_{\rbo}^{-1/2} \leftrightarrow G_{\bfo}^{-2}$ just involves the $v_0$ term, and only when 
there is no resonance. To conclude $PG_{\rbo}^{-\demi}=0$ in both cases.\\ 

We define $G_{\lbo}^{-3/2}$ by symmetry with respect to the $\rbo$ terms as before.
  
\subsection{Resolvent}
Let $G(k) \in \Psi_k^{-2; (-2, 0, 0), \mc{G}}(M, \Omegabht)$ be an operator consistent with all the models we have defined above, and let 
$E(k)=(P+k^2)G(k)-{\rm Id}$. Then  
\[E(k)\in 
\Psi_k^{-\infty,\mc{E}}(M; \Omegabht)\]
where $\mc{E}$ is an index family with $\mc{E}_{\bfo} \geq 1$, $\mc{E}_{\sca} \geq 1$, $\mc{E}_{\rbo} \geq 1/2$, $\mc{E}_{\lbo} \geq 3/2$ and $\mc{E}_{\zf}\geq 2$ and it also has  empty index sets at the other faces. 
By Corollary 2.11 of Part I, the error term $E(k)$ iterates away under Neumann iteration. 
 
The true resolvent is given by $R(k)=G(k)(1+E(k))^{-1}$ and the composition result \cite[Prop. 2.10]{GH} implies
\begin{theo}\label{resolventn=3}  
If $n=3$, $(M,g)$ is asymptotically Euclidean with one end, then  
\[R(k)\in \Psi_k^{-2,(-2,0,0), \mc{R}}(M; \Omegabht)\]
where  where $\mc{R}_{\zf}= -2+\mc{R}'$, $\mc{R}_{\bfo}= -2+\mc{R}'$, 
$\mc{R}_{\sca}=0$, $\mc{R}_{\lbo} = \mc{R}_{\rbo}= -3/2+\mc{R}'$, 
for some integral  index set $\mc{R}' \geq 0$, the index sets at the other faces being empty. 
Moreover, $R_{\zf}^j = G_{\zf}^j$ for $j = -2, -1$, $R_{\bfo}^{-2} = G_{\bfo}^{-2}$, $R_{\rbo}^{-3/2} = G_{\rbo}^{-3/2}$, while $R_{\zf}^0 = G_{\zf}^0$ up to a symmetric finite rank term with range in $\ker_{L^2} P\oplus \cc\psi$
where $\psi$ is the resonant state. 
\end{theo}

\subsection{Comparison with Jensen-Kato}\label{ss:JK} In \cite{JK}, Jensen-Kato
compute the order $-1$ term of the resolvent (that is, $R_{\zf}^{-1}$ in our notation) for $P = \Delta_{\RR^3}  + V$ where $ \Delta_{\RR^3}$ is the Laplacian on flat Euclidean space $\RR^3$.  Their result is (in operator notation)
\begin{equation}
R_{\zf}^{-1} =   \Pi_0 V G_2 V \Pi_0 + \tilde\psi \otimes \tilde\psi, \quad G_2(z, z') = \frac1{24\pi} |z-z'|^2.
\label{PVGVP}\end{equation}
Here we write $\tilde\psi$ for their canonical zero-resonance and $\Pi_0$ for projection onto the $L^2$ null space of $P$. 
Let us verify that this agrees with our $G_{\zf}^{-1}$ defined in \eqref{Gzf-1}, in the case that $M$ is flat $\RR^3$ and the potential function satisfies our conditions. 

First we show that $\tilde\psi$ agrees with our $\psi$. Assuming, as Jensen-Kato do, that $V$ decays at infinity faster than $x^3$, $\tilde\psi$ is characterized by the condition
$$
\Pi_0 V G_2 V \tilde\psi = 0.
$$
We need to adapt their Lemma 2.6 as follows: 

\begin{lem} Suppose that $u, v$ are both in $x^{\beta} L^2$ for some $\beta > 5/2$ and that $v$ is orthogonal to $1$. Then 
if $G_0(z,z')=(4\pi|z-z'|)^{-1}$ is the Green function of $\Delta_{\rr^3}$, 
$$
\langle  u, G_2 v \rangle = \lim_{R \to \infty} \int_{B(0, R)} G_0 u \cdot G_0 v .
$$
\end{lem}
(In \cite{JK} both $u$ and $v$ are required to be orthogonal to $1$, and the inner product on the right is then well-defined without regularization.) 
Using this lemma we compute for any $w \in \ker_{L^2} P$
\begin{multline*}
0 = \langle \Pi_0 V G_2 V \tilde\psi, w \rangle = \langle  V \tilde\psi, G_2 V w \rangle = \lim_{R \to \infty} \int_{B(0, R)} G_0 V \tilde\psi \cdot G_0 V w  \\ = \lim_{R \to \infty} \int_{B(0, R)} \tilde\psi \cdot  w.
\end{multline*}
This shows that $\tilde\psi$ is the same as our $\psi$. Thus it remains to show that the first term in \eqref{PVGVP} agrees with the first term in \eqref{Gzf-1}. 

So consider the operator $\Pi_0 V G_2 V \Pi_0$. 
This operator can be written
\begin{equation}\label{JK-expr}
\begin{gathered}
\frac1{24 \pi} \sum_{i,j}  \psi_i \otimes \psi_j \Big( \int \psi_i(z) V(z) |z-z'|^2 V(z') \psi(z') \, dz \, dz' \Big) \\
= \frac1{24 \pi} \sum_{i,j}  \psi_i \otimes \psi_j \Big( \int \psi_i(z) V(z) (|z|^2 - 2 z \cdot z' + |z'|^2)  V(z') \psi_j(z') \, dz \, dz' \Big).
\end{gathered}\end{equation}
We change back to writing the eigenfunctions as functions rather than half-densities, and use $x$ for a smooth positive function on $\RR^3$ equal to $|z|^{-1}$ for $|z| > 1$.  In this notation
$\psi_i \sim x^2 \sum_j \alpha_{ij} a^j_1$ (see \eqref{alpha}), and $a^j_1 = 0$ for $j \geq d+1$. 
Also recall that $a^j_1$ is an eigenfunction on the sphere with eigenvalue $2$, hence 
$$
a^j_1 = \sum a^j_{1l} \sqrt{\frac{3}{4\pi}} x z_l
$$
in terms of the orthonormal basis $\sqrt{3/4\pi} x z_l$ of this eigenspace. 

First we claim that the $|z|^2$ and $|z'|^2$ terms in \eqref{JK-expr} vanish. To see this we observe that  $V \psi_j = - \Delta \psi_j$ since $\psi_j$ is a zero-mode of $P$. Also, each $\psi_i = O(x^2)$ and, by Theorem~\ref{reg} is conormal at $x=0$. Hence 
$\int \Delta \psi_j = 0$ since $\nabla \psi_j = O(x^3)$ has sufficiently fast decay at infinity. We are therefore left with the cross term involving $z \cdot z'$ in \eqref{JK-expr}. 

Next we claim that 
if either $i > d$ or $j > d$ then the corresponding $z \cdot z'$ term in  \eqref{JK-expr} vanishes. It suffices to assume that $j > d$. 
Then, again using $V \psi_j = - \Delta \psi_j$,
$$
\int z'_k \Delta \psi_j = \int \Delta (z'_k \psi_j) - 2 \nabla_k \psi_j \text{ vanishes.}
$$

We are left with terms where $i, j \leq d$ and involving the cross term
$z \cdot z'$. We may choose a coordinate system in which $\sum_j \alpha_{ij} a^j_1$ is a multiple of $x z_1$. Then 
\begin{equation}\begin{gathered}
\int z'_k \Delta \psi_j = \lim_{R \to \infty} \int_{B_R(0)}  z'_k \Delta \psi_j \\
= \lim_{R \to \infty} \int_{B_R(0)} \Delta (z'_k \psi_j) - 2 \nabla_k \psi_j \\
= \lim_{R \to \infty} \int_{\partial B_R(0)} -\partial_r \Big( r (\frac{z_k}{r}) \big( \sum_j \alpha_{ij} a^j_1 \sqrt{\frac{3}{4\pi}} ( \frac{z_i}{r}) r^{-2} + O(r^{-3}) \big) \Big) \, r^2 d\omega  \\ - 
\lim_{R \to \infty} \int_{\partial B_R(0)} 2 (\frac{z_k}{r})  \Big(  \big( \sum_j \alpha_{ij} a^j_1 \sqrt{\frac{3}{4\pi}} ( \frac{z_i}{r}) r^{-2} + O(r^{-3})  \big) \Big) \, d\omega \\
= -3\sqrt{\frac{3}{4\pi}} \sum_j \alpha_{ij} a^j_1 \int_{s^2} z_i z_k d\omega  
= -\sqrt{12 \pi} \sum_j \alpha_{ij} a^j_1 \delta_{ik}.
\end{gathered}\end{equation}
We see from this that
$$
\frac1{12\pi} \int \Delta \psi_i (z) z \cdot z' \Delta \psi_j(z') =  \sum_{l,m} \alpha_{il} \alpha_{jm} \langle a^l_1, a^m_1 \rangle_{S^2}.
$$
In this form we see, using \eqref{magic1}, \eqref{Gzf-1} and the line below \eqref{D}, that this agrees with our $G_{\zf}^{-1} = R_{\zf}^{-1}$. 

\subsection{Unboundedness of Riesz transform on $L^2$}\label{unboundedness} We see from the results of this section that the term $R_{\zf}^{-1}$ does not vanish if either there is a zero resonance or a zero mode $\psi_k\sim x^{1/2} |dg_b|^{1/2}$ of $P$ (i.e. if $m=1$ in the notation of Theorem~\ref{Riesz}). It follows that the integral \eqref{p-demi} defining $P_{>}^{-1/2}$ does not converge, leading one to suspect that the Riesz transform is not defined --- even on $L^2$. In fact, a simple direct argument shows that this is the case. Suppose for simplicity that the pure point spectrum of $P$ is a single zero mode of $P$ with asymptotic $\psi_k \sim x^{1/2} |dg_b|^{1/2}$. For the Riesz transform to be bounded on $L^2$ we require that there is a constant $C$ such that 
$$
\langle df , df \rangle \leq C  \langle Pf, f \rangle  
$$
for all $f \in H^2(M)\cap (1-\Pi_0)L^2(M)$ (in this case $1-\Pi_0=\Pi_>$ is the projection on the positive spectrum of $P$).  To simplify let us present the argument when there is a single $L^2$ zero mode $\psi_k$ as above. 
Then (after rotating coordinates) $\psi_k \sim z_1/r^3 + O(r^{-3})$. Let  $\chi$ be a cutoff function supported in $\{1/4<|z|<2\}$, equal to $1$ in $\{1/2<|z|<1\}$,  and let $\phi_R :=  R^{-1} \chi(z/R) z_1/r$. Then, as $R\to\infty$, $\| \phi_R \|_{L^2} \sim R^{1/2}$,  
$\| d\phi_R \|_{L^2} \sim R^{-1/2}$ and  
$\langle \phi_R, \psi_k \rangle \sim A$ for some constant $A\not=0$.
 We have $f_R = \phi_R - \langle \phi_R, \psi_k \rangle \psi_k$ in the range of $(\Id - \Pi_0)$ and 
$$
\langle df_R, df_R \rangle \sim |\langle \phi_R, \psi_k\rangle|^{2}||d\psi_k||^2\sim A^2||d\psi_k||^2\not=0.
$$
However, 
$$
\langle P f_R, f_R \rangle \sim  \langle P \phi_R, f_R \rangle \sim \langle P \phi_R, \phi_R \cjd =O(R^{-1}).
$$
So no such constant $C$ can exist. A similar argument with the zero-resonance (taking $\phi_R = \chi(r/R) r^{-1}$) gives the same result and the general case works exactly the same way, using in addition that the eigenvectors corresponding to negative spetrum are vanishing at all order at $\infty$ by Agmon estimates \cite{Ag}.\\

The same arguments clearly works as well for the conic case investigated in the next section and it shows that there exists 
a sequence $f_R\in\Pi_>(L^2(M))$ and $\alpha(R)\to \infty$ as $R\to \infty$ such that $||df_R||^2\geq \alpha(R)\cjg Pf_R,f_R\cjd$, as long as $1/2<m<3/2$ in dimension $n=3$, $0<m<1$ in dimension $n=4$ and $0\leq m<1/2$ in dimension $n=5$.


\section{Resolvent in dimension 3: asymptotically conic manifolds}\label{sec6}

It is interesting to see what happens in the three-dimensional case when the manifold is asymptotically conic, rather than asymptotically Euclidean. In this case there may be models at $\zf$ at non-integral orders between $-2$ and $0$, i.e. nontrivial $R_{\zf}^\alpha$ for $-2 < \alpha < 0$, $\alpha \neq -1$. Indeed, the fact that there is a nonzero term $R_{\zf}^{-1}$ at the particular order $-1$ in the previous section is related to the precise value of the second eigenfunction of the Laplacian on $S^2$ with the standard metric, and is not at all typical in the class of three-dimensional asymptotically conic manifolds. 

Rather than attempt a comprehensive treatment, we consider only a couple of very special cases. 
We give very few details and only indicate the most interesting features of these examples. 

\subsection{One zero resonance and no zero modes}\label{zr-conic-3}
We assume that $(M,g)$ is asymptotically conic and recall the set 
$N_{\pl M}$ from \eqref{npl}. Let us assume that  $P$ has 
one unique resonant state $\psi$ at zero energy with asymptotic  $\psi(x,y)\sim a_\nu(y) x^{\nu-1} |dg_b|^{\demi}$
for some $\nu\in [1/2,1)\cap N_{\pl M}$, $\eps>0$ and $a_\nu\in E_{\nu}$.
In terms of the operator $P_b$, this means that $P_b$ has a unique $L_b^2$ (normalized) zero-mode $\varphi = x \psi$ with asymptotic 
\begin{equation}\label{varphiasympt}
\varphi(x,y)=a_\nu x^{\nu}+\sum_{\nu<\nu_i\leq \nu+2}x^{\nu_i}a_{\nu_i}(y)+O(x^{\nu+2+\eps})\end{equation}
where $a_{\nu_i}\in E_{\nu_i}$.

We choose $G_{\zf}^0$ as in Section~\ref{dim3}:
\[G_{\zf}^0=(xx')^{-1}\Big(Q_b+u\otimes \varphi+\varphi\otimes u\Big)\]
where $Q_b$ is the generalized inverse of $P_b$, i.e. $P_bQ_b=\Id-\varphi\otimes\varphi$ and $u\in x^{-\nu-\eps}L^2_b$ is 
a function such that $P_bu=\varphi$. The existence of such a $u$ with asymptotic  
\begin{equation}\label{asymptu}
u=b^-_{\nu}x^{-\nu}+\sum_{1/2\leq \nu_i<\nu}(b^-_{\nu_i}x^{-\nu_i}+b^+_{\nu_i}x^{\nu_i})-
\frac{a_{\nu}}{2\nu}x^{\nu}\log x +\sum_{\nu<\nu_i\leq 1}b^+_{\nu_i}x^{\nu_i}+O(x^{1+\eps}),\end{equation}
where the $\nu_i$ runs over the set $N_{\pl M}$ and $b^\pm_{\nu_i}\in E_{\nu_i}$, is ensured by Theorem \ref{rit} and Theorem \ref{reg}. 

We define $G_{\bfo}^{-2}$ similarly to the Euclidean case with a resonance, that is 
$G_{\bfo}^{-2}=\kappa\kappa' Q_{\bfo}$ where 

\begin{multline*}
Q_{\bfo} = \Bigg[\sum_{j=0}^\infty \Pi_{E_j} \Big(I_{\nu_j}(\kappa)K_{\nu_j}(\kappa')H(\kappa'-\kappa)+
I_{\nu_j}(\kappa')K_{\nu_j}(\kappa)H(\kappa-\kappa')\Big)\\
+ c_{\nu}a_\nu(y)a_{\nu}(y')K_{\nu}(\kappa)K_{\nu}(\kappa') \Bigg]  
\left|\frac{d\kappa dyd\kappa'dy'}{\kappa\kappa'}\right|^\demi\end{multline*}
where $c_\nu\in\cc$ is to be determined to get matching with $G_{\zf}^0$.
Considering the asymptotic expansion of $G_{\bfo}^{-2}$ at $\zf$ gives
\begin{multline}\label{gbfonu}
G_{\bfo}^{-2}=\Gamma(\nu)^22^{2\nu-2}(\kappa \kappa')^{1-\nu}a_{\nu}(y)a_{\nu}(y')+
(\kappa\kappa')^{-1}I_{\rm ff}(Q_b)\\
+c_\nu a_{\nu}(y)a_{\nu}(y')\frac{\Gamma(\nu)\Gamma(-\nu)}{4}
( \kappa^{1-\nu}{\kappa'}^{1+\nu} + \kappa^{1+\nu}{\kappa'}^{1-\nu} )+O(\rho_{\zf}^{2+\eps}).
\end{multline}
The $c_\nu$ term has to match with the $(xx')^{-1}(u\otimes\varphi+\varphi\otimes u)$ term of $G_{\zf}^0$. A short computation proves that they match as long as $c_\nu=-2/(\Gamma(\nu)\Gamma(1-\nu)||a_\nu||^2)$. This will force 
a term at order $-2\nu$ to be added at $\zf$. We set
\[G_{\zf}^{-2\nu}:=\psi\otimes \psi\]
where $\psi=2^{\nu-\demi}\sqrt{\Gamma(\nu)/\Gamma(1-\nu)}||a_{\nu}||^{-1}x^{-1}\varphi$ is chosen (up to sign) to match with the $(\kappa \kappa')^{1-\nu}$ term of (\ref{gbfonu}). In other words, $\psi$ is the unique function in the null space of $P$
with asymptotic 
\[\psi=2^{\nu-\demi}\sqrt{\frac{\Gamma(\nu)}{\Gamma(1-\nu)}}x^{\nu}\phi_\nu(y)+O(x^{\nu+\eps})\]
where $\phi_\nu$ is a function in $E_{\nu}$ with $L^2(\pl M)$ norm $1$. Note that the order of this term is  $-2\nu \in [-1, 2)$ under our assumptions on $\nu$. 

These models can be completed to a parametrix $G(k)$ with error term that iterates away, showing that the resolvent is phg on $\MMksc$ and that $R_{\zf}^{-2\nu} = G_{\zf}^{-2\nu} = \psi\otimes \psi$. 

\subsection{One zero mode and no zero resonance} 
The next example we consider is that  $P$ has 
one zero mode $\psi(x,y)\sim a_\nu(y)x^{\nu-1}$ for some $\nu\in(1,3/2]$, $a_\nu\in E_\nu$.
In terms of the operator $P_b$, this means that $P_b$ has a unique (up to sign) $L^2$-normalized zero mode $\varphi = x \psi$ with asymptotic 
\[
\varphi(x,y)=a_\nu x^{\nu}+\sum_{\nu<\nu_i\leq \nu+2}x^{\nu_i}a_{\nu_i}(y)+O(x^{\nu+2+\eps})\]
where $a_{\nu_i}\in E_{\nu_i}$.
We may take $\psi$ to be $L^2$-normalized and write $\psi=cx^{-1}\varphi$ for $c:=||x^{-1}\varphi||^{-1}_{L^2}$.

We follow the method of Section \ref{dim3}. We set the first term at $\zf$ to  
\[G_{\zf}^0:=(xx')^{-1}\Big(Q_b +\chi\otimes \varphi +\varphi\otimes \chi\Big)\]
where $Q_b$ is the generalized inverse of $P_b$ and $\chi\in x^{\nu-2-\eps}L^2_b$ is a function such that 
$P_b\chi=x^{-1}\psi^\perp$ where $\psi^\perp:=x\varphi-\cjg x\varphi,\psi\cjd\psi=x\varphi-c\psi$ (recall 
$\psi=cx^{-1}\varphi$ is the $L^2$ normalized eigenvector). We thus get $PG_{\zf}^0=\Id-\psi\otimes\psi$ 
and we set 
\[G_{\zf}^{-2}:=\psi\otimes \psi, \quad G_{\bfo}^{-2}:=\kappa\kappa'Q_{\bfo}\] 
where $Q_{\bfo}$ is given in (\ref{bfo-conic-model}). 

The $G_{\zf}^{-2}$ term requires a matching term at order $\nu - 3$ at $\rbo$ which we take to be 
$$
G_{\rbo}^{\nu - 3} = \frac{\kappa'K_{\nu}(\kappa')}{\Gamma(\nu)2^{\nu-1}}ca_{\nu}(y')\psi(z).
$$
This in turn requires a term at order $2\nu - 2$ at $\zf$, from the $\kappa^\nu$ term in $K_{\nu}(\kappa)$ as $\kappa \to 0$. To match this with $G_{\bfo}^{\nu - 3}$, we note that 
there exists $u\in x^{-\nu-\eps}L^2_b$ such that
$P_bu=\varphi$ and with asymptotic 
\[u(x,y)= a_{\nu}(y)x^{-\nu}\gamma^{-1}+\sum_{1/2\leq \nu_i\leq 1}\gamma^+_{\nu_i}(y)x^{\nu_i}+O(x^{1+\eps}) 
\] with $\gamma:=-2\nu||a_\nu||_{L^2(\pl M)}^2$ (we have used Proposition \ref{prop:complementary} so that the $x^{-\nu_i}$ powers do not show up for $1/2\leq \nu_i <\nu$).
Thus we have $P(x^{-1}u)=x\varphi$ and $P(x^{-1}(\chi-u))=\psi^\perp-x\varphi=-c\psi$.
Then we define 
\[G_{\zf}^{2\nu-2}:=-\gamma c2^{-2\nu}\frac{\Gamma(-\nu)}{\Gamma(\nu)}\Big(x^{-1}(\chi-u))\otimes\psi +\psi\otimes(x^{-1}(\chi-u))\Big)\]
\[G_{\zf}^{2\nu-4}:=-\gamma c^2\frac{\Gamma(-\nu)}{\Gamma(\nu)}2^{-2\nu}(\psi\otimes\psi)\]
so that $PG_{\zf}^{2\nu-2}=-G_{\zf}^{2\nu-4}$ and they match
with the ${G}^{\nu-3}_{\rbo}$ term at $\rbo\cap \zf$. 

Following the same method used for the asymptotically Euclidean case above, we can complete this to a parametrix so that the error term iterates away, and we find that the resolvent is again polyhomogeneous, with
$R_{\zf}^j = G_{\zf}^j$ for $j = -2, 2\nu - 4$. 

\subsection{A limiting case} We have not discussed what happens in the case $\nu = 1$. This corresponds to a zero-resonance that fails only logarithmically to be square-integrable. In Section~\ref{zr-conic-3} we saw that for $\nu \in [1/2, 1)$ there is a nonzero term $R_{\zf}^{-2\nu}$, but with a coefficient $\Gamma(1-\nu)^{-1}$ that vanishes at $\nu =1$, so it is not clear what to expect at $\nu = 1$. (We can also note that the simple bound $\| R(k) \|_{L^2 \to L^2} \leq k^{-2}$ valid near $k=0$ prevents any behaviour more singular than $k^{-2}$ in the kernel of $R(k)$ as $k \to 0$, while a term at order $k^{-2}$ can only be expected if there is an atom in the spectral measure at $k=0$.)

What happens is that the resolvent becomes non-polyhomogeneous in this case and there is a term of order $k^{-2} (\log k)^{-1}$ in the expansion of the resolvent at $\zf$. Even though the resolvent in this case fails to be non-polyhomogeneous, the polyhomogeneous calculus can be used to construct a weak parametrix with a $O(1/\log k)$ error term. 
Rather than do this here, instead we move to the four-dimensional case since the same phenomenon occurs there in the asymptotically \emph{Euclidean} situation, which is perhaps of greater interest and importance. This is the topic of the following section.

\section{Resolvent in dimension $4$ with a zero-resonance - asymptotically Euclidean case}
Assume now that $(M,g)$ is four-dimensional and asymptotically Euclidean. 
To simplify the exposition, we assume \eqref{assump1}  as before. 
As mentioned above, this case is a special case in which the resonance only fails logarithmically to be square-integrable, and this leads to non-polyhomogeneous behaviour of the resolvent. 
Our strategy is to find a polyhomogeneous kernel $G(k)$ on $M^2_{k,\sca}$ such that
\[(P+k^2)G(k)=\log k \cdot \Id + E(k)\]
for some polyhomogeneous function $E(k)$ which is the integral kernel of a compact operator uniformly down to $k=0$, with 
$L^2 \to  L^2$ operator norm uniformly bounded. This means that $\Id + E(k) /\log k$ can be inverted using the Neumann series for small $k$, and this yields the resolvent $(G(k)/\log k) (\Id + E(k)/\log k)^{-1}$ which we see has an expansion in powers of $1/\log k$ (in the interior of $\zf$, say).

To avoid problems coming from the nature of the function $\log k$ near $\bfo\cap\zf$, 
we will now switch to use the boundary defining function $\kappa'=k/x'$
for $\zf$, and $x'$ for $\bfo$ (at least outside $\bfo\cap\rbo$) and for $\rbo$, both of which commute with $P$.  
Following our previous notations, we shall denote by $G_{\rm f}^{i,j}$ the term defined on the face $\rm f$ which corresponds
to the $\rho_{\rm f}^i\log(\rho_{\rm f})^j$ coefficient in the asymptotic expansion at the face ${\rm f}$, where 
$\rho_{\rm f}$ is the previously fixed boundary defining function of the interior of the face $\rm f$.

\subsection{Diagonal and scattering face terms}
The diagonal singularity of $G(k)$ is encoded in a symbolic term which is the symbolic inverse of $\log(k)\sigma(P+k^2)$ and
the leading term at the scattering face $\sca$ is the inverse of $\log(k)(\Delta+k^2)$ where $\Delta$ is the flat 
laplacian on the fibers of $\sca$ induced by $\Delta_g$.

\subsection{Terms $G_{\rm zf}^{0,1}$, $G_{\rm zf}^0$}

A resonant state $\psi$ correspond to an $L^2_b$ normalized element $\varphi\in\ker P_b$ in $L^2_b(M)\setminus x^{\eps}L^2_b(M)$ ($\eps>0$ small) by the relation $\psi\in \cc x^{-1}\varphi$, it has asymptotic 
$\varphi=ax+O(x^2)$ for some $a\in\cc$.
We follow the method applied for $n=3$ and thus we will only give few details. 
There is a generalized inverse $Q_b$ for $P_b$ such that
\[P_bQ_b={\rm Id}-\varphi\otimes\varphi.\]
By Theorem~\ref{rit}, there exists $u\in x^{-1-\eps}L^2_b$ such that 
$P_bu=\varphi$  and with asymptotics
\[u=bx^{-1}-\frac{a}{2}x\log x +O(x^2) \]
for some $b,c\in\cc$ satisfying $ab=-(2{\rm vol}(S^3))^{-1}$ from Green's formula applied to 
$$
1 = \lim_{\epsilon \to 0} \int_{M \cap \{ x \geq \epsilon \}} \langle P_b u, \varphi \rangle - \langle u, P_b \varphi \rangle.
$$
We set 
\[G_{\zf}^{0,1}=(x x')^{-1}  (Q_b+u\otimes \varphi+\varphi\otimes u)\]
so that $PG_{\zf}^{0,1}= \Id$ on $\zf$. Note that $I_{\rm ff}(Q_b)$, the restriction of $Q_b$ to $\ff \equiv \zf \cap \bfo$,  is given by (\ref{nkernelff}) with $n=4$. 
Like the case $n=3$ in (\ref{asymptilQ}), we have the asymptotic at the right boundary 
\begin{equation}\label{asymptrightb}
(xx')^{-1}(Q_b+u\otimes\varphi+\varphi\otimes u)=b{x'}^{-2}\varphi(z)+\beta\varphi(z)+O(x')\end{equation}
for some $\beta\in\cc$.
At $\bfo\cap\zf$, we have 
\begin{equation}\label{gzfn=4reso}
G_{\zf}^{0,1}=(x x')^{-1}  I_{\rm ff}(Q_b)-\frac{x^{-2}+{x'}^{-2}}{2{\rm Vol}(S^3)}+O(\rho^{-1}_{\bfo}\log^2(\rho_{\bfo})).
\end{equation}
We shall also define (the tilde indicates that it will be corrected below) 
\[\til{G}_{\zf}^0=\log(x')(x x')^{-1}(Q_b+u\otimes\varphi+\varphi\otimes u).\]
which satisfies $P\til{G}_{\zf}^0=\log (x'){\rm Id}.$
Note that $P(\log(\kappa')G_{\zf}^{0,1}+G_{\zf}^0)=\log(k){\rm Id}$ since $\log\kappa'+\log x'=\log k$, moreover
the same argument shows that
$G_{\zf}^{0,1}, G_{\zf}^0$ match with the diagonal singularity term.

\subsection{Term $G_{{\rm bf}_0}^{-2,1}$ and $G_{{\rm bf}_0}^{-2}$}
With this new boundary defining function for $\bfo$, we write $P+k^2={x'}^2\kappa P_{\bfo}\kappa^{-1}$
and thus we set $G_{\bfo}^{-2,1}=\frac{\kappa}{\kappa'}Q_{\bfo}$ where
\[Q_{\bfo}=\sum_{j=0}^\infty \Pi_{E_j}(y,y') \Big(I_{1+j}(\kappa)K_{1+j}(\kappa')H(\kappa'-\kappa)+
I_{1+j}(\kappa')K_{1+j}(\kappa)H(\kappa-\kappa')\Big).\]
The term $\log(x'){x'}^{-2}G_{\zf}^{-2,1}$ matches with the $\log(x)(xx')^{-1}Q_b$ part of $\til{G}_{\zf}^{0}$
and certainly match with $\log(\kappa')G_{\zf}^{0,1}$ at $\zf\cap\bfo$.
We can not add a finite rank term $K_{1}(\kappa)K_1(\kappa')$ to $G_{\bfo}^{-2,1}$ like we did 
for the case $n=3$ to match $G_{\zf}^{0}$ since this would imply a term of order $\rho^{-2,1}_{\zf}$ 
at $\zf$ to be added, and then 
we would not be able to end the construction of a parametrix with a sufficient error.
However we can instead put a term of order $\rho_{\bfo}^{-2}$ at $\bfo$ that we call
$G_{\bfo}^{-2}$, we define it by 
\[G_{\bfo}^{-2}:=\log(\kappa')\frac{\kappa}{\kappa'}Q_{\bfo}
+c_0 \frac{\kappa}{\kappa'} K_{1}(\kappa)K_1(\kappa')
\]
for some $c_0\in\cc$ and we get $I_{\bfo}({x'}^{-2}(P+k^2))G_{\bfo}^{-2}=\log(\kappa'){\rm Id}$.
Using the asymptotic
\[K_1(z)=z^{-1}+\frac{1}{2} z\log z+\alpha z+O(z^3\log z), \quad \alpha:=-\frac{2\log(2)+1-2\gamma}{4}\] 
at $z=0$ with $\gamma$ Euler's constant, we find 
\begin{equation}\label{gbfo2-1}
G_{\bfo}^{-2}=c_0{\kappa'}^{-2} + \log(\kappa')\Big(\frac{c_0}{2}(t^2 +1)+tQ_{\bfo}\Big) 
+ c_0\Big(\alpha(1+t^2)+\frac{t^2}{2}\log(t)\Big)+
O(\rho_{\zf})\end{equation}
where $t:=\kappa/\kappa'=x'/x$ is smooth and non vanishing outside $\rbo,\lbo$ but near $\zf\cap\bfo$.  The term of order ${x'}^{-2}\log(\kappa')$ matches with the ${\rm vol}(S^3)$ term of $\log(\kappa')G_{\zf}^{0,1}$
if and only if $c_0=-1/{\rm vol}(S^3)=2ab$.
The consistency of these two terms with the diagonal singularity is again straightforward.

\subsection{Term $G_{\rm zf}^{-2}$}
The most singular term of the asymptotic (\ref{gbfo2-1}) at $\zf$ forces to add a term of order 
${\kappa'}^{-2}$ at $\zf$, we call it $G_{\zf}^{-2}$ and define it by
\[G_{\zf}^{-2}=-\frac{\psi\otimes\psi}{{x'}^2}, \quad \textrm{ with } \psi\in\cc(x^{-1}\varphi), \textrm{ such that } 
\psi\sim\sqrt{\frac{1}{{\rm vol}(S^3)}}.\]
(Recall that now we are using $\kappa'$, not $k$, as the boundary defining function for $\zf$.) 
Note that ${\kappa'}^{-2}G_{\zf}^{-2}$ certainly matches with ${\kappa'}^{-2}\log(\kappa')G_{\bfo}^{-2,1}$ and ${x'}^{-2}G_{\bfo}^{-2}$; moreover it solves $PG_{\zf}^{-2}=0$ at $\zf$.

\subsection{Term $G_{\rm zf}^{0}$ corrected}
We need to correct $G_{\zf}^0$ so that it matches with the $\bfo$ terms. We just set 
\[G_{\zf}^{0}:=\log(x')(x x')^{-1}(Q_b+u\otimes\varphi+\varphi\otimes u)+(\theta\otimes \psi+\psi\otimes \theta)\] 
where $\theta$ is chosen to be a polyhomogeneous function satisfying 
\[\theta = \frac{x^{-2}\log x}{2\sqrt{{\rm vol}(S^3)}}-\sqrt{\frac{1}{{\rm vol}(S^3)}}\alpha x^{-2}+O(x).\]
A short computation shows that $G_{\zf}^0$ matches with the $\bfo$ terms and $PG_{\zf}^0=\log (x'){\rm Id}+(P\theta)\otimes\psi$.

\begin{remark}\label{psifp}
In principle, to improve the parametrix and get a better error term at $\zf$, it would be necessary to have
$P\theta=\psi$ so that $PG_{\zf}^{0}=(\log x'){\rm Id}-k^2{\kappa'}^{-2}G_{\zf}^{-2}$. Using the Theorem \ref{rit}, we can find a function $v\in x^{-1-\eps}L^2_b$ such that $P_bv=x^{-1}\psi$, or equivalently $P(x^{-1}v)=\psi$ and one can check  
that $v$ satisfies the asymptotic $v=\frac{1}{2\sqrt{{\rm Vol}(S^3)}}x^{-1}\log x +\mu x^{-1}+O(x\log x)$ for some $\mu\in\cc$. However we are not able to prove that $\mu=-\alpha/\sqrt{{\rm Vol}(S^3)}$ in order to take $\theta:=x^{-1}v$, and we believe that 
this actually does not hold in general. Notice in addition that Green's formula
applied to $\int_{x>\eps} P(x^{-1}v)\psi$ shows that 
\begin{equation}
||\psi||^2_{\rm fp} :=\textrm{fp}_{\eps\to 0} \int_{x>\eps} |\psi|^2=\demi-\mu\sqrt{{\rm Vol}(S^3)}
\label{fp}\end{equation}
where fp means finite part (i.e. we keep the $\eps^0$ coefficient in the asymptotic expansion at $\eps \to 0$). 
\end{remark}

\subsection{Terms at ${\rm rb}_0$}
We first set 
\[G_{\rbo}^{-2}:=-{\kappa'}^{-1}K_1(\kappa')x^{-1}\varphi(z)\]
which matches with the 3 $\zf$ terms using (\ref{asymptrightb}). 
There is a term at order $-1$ at $\rbo$ coming from the second asymptotic term
in the expansion of $\psi$ in $G_{\zf}^{-2}$: writing 
\[\psi=1/\sqrt{{\rm Vol}(S^3)}+\phi_1(y)x+O(x^2\log x)\]
for some $\phi_1\in E_1$ (recall $E_j=\ker(\Delta_{S^3}-3)$), we set
\[G_{\rbo}^{-1}:=\demi K_2(\kappa')\phi_1(y')\psi(z)\]
Then using again (\ref{asymptrightb}) we see that the term
\[G_{\rbo}^{0,1}:=\beta x^{-1}\varphi(z)\kappa'K_1(\kappa')\]
matches with the $\zf$ terms.
Notice that the $3$ terms at $\rbo$ solve $PG_{\rbo}^*=0$ and are consistent with the $\bfo$ terms.

\subsection{Term at $\lbo$}
Here we may use again the function $k$ as a defining function for $\lbo$, then we set
\[G_{\lbo}^{-2}:=-{\kappa}K_1(\kappa){x'}^{-1}\varphi(z').\]
Notice that $G_{\lbo}^{-2}$ solves the model equations at the $\lbo$ face. 

\subsection{Resolvent}
One can construct an operator $$G(k)\in \log k \cdot \Psi_k^{-2; (-2, -2, 0), \mc{G}}(M, \Omegabht) + \Psi_k^{-2; (-2, -2, 0), \mc{G}}(M, \Omegabht)$$ consistent with all the models we have defined above, where  $\mc{G}$ is an index set such that
\[\mc{G}_{\sca}=0, \quad \mc{G}_{\bfo}\subset (-2,0)+\mc{G}', \quad \mc{G}_{\zf}\subset-2+\mc{G}',\quad \mc{G}_{\rbo}=\mc{G}_{\lbo}=-2+\mc{G}'\] 
where $\mc{G}'>0$ is some integral index set (including log terms). 
The error defined by 
$E(k)=-(P+k^2)G(k)+(\log k){\rm Id}$ is in $\Psi_k^{-\infty,\mc{E}}$ where $\mc{E}$ is an index set (integral 
except log terms) such that for any $\eps>0$
\begin{equation}\label{ek}
\mc{E_{\zf}}\geq 0, \quad \mc{E}_{\bfo}\geq 1-\eps, \quad \mc{E}_{\rbo}\geq 0, \quad \mc{E}_{\lbo}\geq 1, \quad \mc{E}_{\sca} \geq 1
\end{equation}
with the empty set at the other faces.  

We can now state a result concerning the expansion of the resolvent as $k \to 0$ in two different senses:

\begin{theo}\label{n4resonance}
Assume that $n=4$, that $M$ is asymptotically Euclidean, and that  $\ker_{L^2} P=0$ but $P$ has a zero-resonance $\psi$ at $0$ normalized as above. Let $R(k)$ denote the resolvent $(P + k^2)^{-1}$ of $P$. 

(i) The kernel $k^2 R(k)$ has an expansion 
\begin{equation}
k^2 R(k) = \sum_{j=0}^\infty (\log k)^{-j} R_j(k), 
\label{R-series-1}\end{equation}
with $R_j(k)$ in the calculus: 
\begin{equation}\begin{gathered}
R_j(k) \in \Psi_k^{-2; (0,0,0), \mc{R}^j}(M, \Omegabht), \\
\mc{R}^j_{\sca}=\nn_0,\quad \mc{R}^j_{\zf}\subset (0,0) \cup {\mc{R}'}^j, \quad \mc{R}^j_{\rbo},\mc{R}^j_{\lbo}\subset{\mc{R}'}^j,\quad
\mc{R}^j_{\bfo}=(0,0) \cup{\mc{R}'}^j
\end{gathered}\end{equation}
for some index set ${\mc{R}'}^j>1-\eps$ for any $\eps>0$. The terms in the series are eventually bounded operators on $L^2$, uniformly in $k$, and the series converges in the norm $L^\infty([0, k_0]; \| \cdot \|_{L^2 \to L^2})$ for sufficiently small $k_0 > 0$. 

(ii) In addition, for any compact set $K$ of $\{(z,z')\in M^\circ\x M^\circ, z\not=z'\}$, and for any $N\in\nn,\eps>0$ the kernel $k^2 R(k;z,z')$ of $k^2 R(k)$ has an expansion of the form   
\begin{equation}
k^2 R(k;z,z')= \sum_{\ell=0}^N\sum_{j=-J(\ell)}^\infty R_{\ell,j}(z,z')k^\ell (\log k)^{-j} +O(k^N(\log k)^{J(N+1)}) 
, \quad (z,z')\in K
\label{R-series-2}\end{equation}
for some $J:\nn\to \nn$ and $R_{\ell,j}\in C^\infty(K)$, which converges in $L^\infty([0, k_0]; C^\infty(K))$. 
The kernels $R_{0,j}$ satisfy $R_{0,j}=0$ for $j\leq 0$ and 
\[R_{0,j}(z,z')= (-1)^{j} \omega^{j-1}\psi(z)\psi(z')
\]
for $j>0$, where the real number $\omega$ is given by \eqref{omega}.

\end{theo}

\begin{proof}
Let us define $E'(k):=E(k)/ \log k$.
We first prove that, for $j\geq n/2+1$, $E'(k)^{j}$ is a uniformly bounded (down to $k=0$) family of  Hilbert-Schmidt operators with norm converging to $0$ as $k\to 0$.  In view of (\ref{ek}) and the composition theorem (Proposition 2.10  of Part I \cite{GH}), the kernel $E_j(k;z,z')$ of $E(k)^j$ is bounded by
\begin{equation}
|E_j(k;z,z')|\leq Cx^{1-\eps}\Big(\frac{x'}{x'+k}\Big) |dg_b dg_b'|^{1/2}
\label{Ej-est}\end{equation}
where $C$ is uniform in $k\in(0,1)$: indeed $E(k)^j$ is an operator in $\Psi_k^{-\infty,\mc{E}^j}$
for some index set $\mc{E}^j$ satisfying same bounds as in (\ref{ek}) but with $\mc{E}_{\sca}^j>n/2$ in addition. (See Definition 2.7 of Part I for why we need the vanishing factor of order $ > n/2$ at $\sca$.) 
We consider the $L^2$ norm of the kernel of $E_j(k,z,z')$ on $M^2$: let us fix $k$, then in the region 
$k\geq x'$ we have (it suffices to work in $x,x'\in(0,1)$)
\[ \int_{S^3}\int_{S^3}\int_{0}^1\int_{0}^k \Big| E_j(k;z,z') \Big|^2\frac{dx'}{x'}\frac{dx}{x}dydy' \leq C'\int_{0}^1\Big(\frac{x'}{x'+1}\Big)^2\frac{dx'}{x'}\leq C''\] 
where $C',C''$ are constants. It remains to treat the part $x'>k$. Then we remark that $|\log k|^{-1}<|\log x'|^{-1}$
and 
\[|E_j(k;z,z')(\log k)^{-1}|\leq C|\log k|^{-1/3}x^{1-\eps}|\log x'|^{-2/3}\]
which is in $L^2_b(M\x M)$ with $L^2$ norm bounded by $|\log k|^{-1/3}$ (this is because $|\log x'|^{-1-\epsilon}\in L^1((0,1/2),dx'/x')$ for any $\epsilon > 0$). By iterating this, it follows that $\| E'(k)^j \|_{L^2 \to L^2}$ is bounded by $A C^j (\log k)^{-j/3}$ for some $A > 0$ and for all $j > n$, where $C$ is independent of $j$.  

We can apply similar reasoning to the kernel $k^2 G(k)/\log k  \circ E'(k)^j$ for a fixed $j \geq n/2+2$. Write $G(k) = G_1(k) + G_2(k),$ with $G_i(k) \in \Psi_k^{-2; (-2, 0, 0), \mc{G}}(M; \Omegabht)$. Then using the composition formulae, $k^2 G_i(k) \circ E(k)^j$ is in the calculus with index sets $\geq 0$ at $\zf$, $\bfo, \lbo, \rbo$ and $\geq n/2+2$ at $\sca$. Instead of \eqref{Ej-est} we can bound the kernel by 
$$
\Big(\frac{x'}{x'+k}\Big)\Big(\frac{x}{x+k}\Big) |dg_b dg_b'|^{1/2}
$$
(since the kernel, as a multiple of $|dg_b dg'_b|^{1/2}$, vanishes to at least second order at $\sca$, and to any order at $\lb, \rb, \bfc$). 
Using similar reasoning we deduce a bound $|\log k|^{4/3}$ on the $L^2$ norm of this kernel. But we have an extra factor $(\log k)^{-j}$ in $k^2 G(k)/\log k  \circ E'(k)^j$ which yields a  uniform bounds for $\| k^2 G(k)/\log k \circ  E'(k)^{j} \|_{L^2 \to L^2}$ for any $j \geq n/2+2$. Combining these results we see that the series \eqref{R-series-1} is convergent in operator norm, proving the first part of the theorem.

To prove the second part, we use the fact that the $L^\infty$ norm of the kernel of an operator coincides with the $L^1 \to L^\infty$ operator norm. With a similar argument to that above, we deduce that  $\| E'(k)^j \|_{L^2 \to L^\infty}$ and $\| k^2 G(k)/\log k \circ  E'(k)^j \|_{L^1 \to L^2}$ are bounded uniformly in $k$ for any $j > n/2$. (These are strictly weaker results since we only need the $L^2 L^\infty$ norm on the kernel instead of the $L^2 L^2$ norm, and we already know from the calculus that these operators have bounded kernels.) Composing these inequalities with $\| E'(k)^j \|_{L^2 \to L^2} \leq A C^j (\log k)^{-j/3}$, we conclude that the $L^\infty$ norm of the kernel of  $k^2 G(k) E'(k)^j/\log k$ is bounded by $A C^j (\log k)^{-j/3}$ for large $j$ and therefore the series \eqref{R-series-2} converges in $L^\infty$ for $k$ small. Moreover, the same estimates hold if we apply arbitrary smooth differential operators on $K$ to the kernel of  $G(k) E'(k)^j$, so we get convergence of \eqref{R-series-2} in $L^\infty([0, k_0]; C^s(K))$ for all $s$.

To compute the terms of order $k^{-2}(\log k)^{-j}$ in the expansion of $R(k)$ at $\zf$, we notice that, at $\zf$, 
\[G(k)=\frac{G_{\zf}^{-2}}{k^2 \log k}+O(1), \quad E(k)^j=(E_{\zf}^0+O(k))^j=(E_{\zf}^0)^j+O(k)\] 
where $E_{\zf}^0:=PG_{\zf}^{0}+k^2G_{\zf}^{-2}=(P\theta-\psi)\otimes \psi$ is the term of order $k^0$
in the expansion of $E(k)$ at $\zf$.
Thus the composition $G_{\zf}^{-2}(E_{\zf}^0)^j$ is the coefficient of order $k^{-2}(\log k)^{-j-1}$
of $R(k)$ at $\zf$, and it is exactly 
\[-\big( \langle \psi,P\theta-\psi \rangle \big)^{j} \psi\otimes\psi.\] 
Now the constant $\omega:=\cjg \psi,P\theta-\psi\cjd$ can be computed through Green's formula 
on balls $\{x>\eps\}$ and taking the limit $\eps\to 0$, this gives 
\begin{equation}
\omega = \demi+\alpha -||\psi||_{\rm fp}^2=-\frac{2\log(2)-1-\gamma}{4}-||\psi||^2_{\rm fp}.
\label{omega}\end{equation}
This completes the proof of the theorem. 
\end{proof}

Note that, using Remark \ref{psifp}, this constant $\omega$ vanishes if and only if $\theta-x^{-1}v\in \cc\psi$ where
$v$ is the unique function (modulo $\cc (x^{-1}\psi)$) such that $P_bv=x^{-1}\psi$.


\section{Riesz transform}
In this section we prove Theorem~\ref{Riesz}. We begin with a preparatory lemma. 

\subsection{A computation}

Let us define the function for $\nu> 1/2$
\[F(\kappa):={\kappa}^{\nu}K_{\nu}(\kappa)\]
which is smooth on $(0,\infty)$, with $F(0)=-2^{\nu-1}\Gamma(\nu)$ 
and $F(\kappa)=F(0)+O(\kappa^2)+O(\kappa^{2\nu})$ as $\kappa\to 0$. 
\begin{lem}\label{comp}
We have the identity for $\nu>1$
\[\int_{0}^\infty \kappa^{-2}(F(\kappa)-F(0))d\kappa=\frac{2^\nu\pi^\demi\Gamma(\nu+\demi)}{-2\nu+1}\]
\end{lem}
\begin{proof} We first have by integration by parts for $\eps>0$ small
\begin{equation}\label{integral}
\int_{\eps}^\infty \kappa^{-2}(F(\kappa)-F(0))d\kappa=\eps^{-1}(F(\eps)-F(0))+\int_{\eps}^\infty \kappa^{-1}\pl_\kappa
F(\kappa)d\kappa\end{equation}
but since $z^{\nu}K_{\nu}(z)=c_\nu\mc{F}_{t\to z}((1+t^2)^{-\nu-\demi})$ with $c_\nu:=\Gamma(\nu+\demi)2^{\nu-1}\pi^{-\demi}$ we get by another integration by parts
\[\frac{\pl_\kappa F(\kappa)}{\kappa}=-ic_{\nu}\int_{\rr}\frac{e^{-i\kappa t}}{\kappa}t(1+t^2)^{-\nu-\demi}dt=
\frac{c_{\nu}}{-2\nu+1}\int_\rr e^{-i\kappa t}(1+t^2)^{-\nu+\demi}dt.\]
Using that formula and making $\eps\to 0$ in (\ref{integral}) yields
\[\int_{0}^\infty \kappa^{-2}(F(\kappa)-F(0))d\kappa=\pl_\kappa F(0)+\frac{2\pi c_{\nu}}{-2\nu+1}\]
which proves the lemma since $\pl_kF(0)=0$. 
\end{proof}

\subsection{Proof of Theorem~\ref{Riesz}}
In this section, we follow essentially the Section 5.2 of Part I \cite{GH}. 
First recall from the Introduction that the Riesz transform $T$ by 
$T = d \circ (P_>)^{-1/2}$, or equivalently (when it is defined on $L^2$)
\begin{equation}
T = \frac{2}{\pi} d \int_0^\infty R(k) \Pi_> \, dk
\label{T}\end{equation}
where $\Pi_>$ is the spectral projection onto $(0, \infty)$ for the operator $P$. We shall also use  $\Pi_0$ to denote the projection onto the zero eigenspace. 
By the arguments of Part I, Sec. 5.2, we just have to study $L^p$ boundedness of 
\begin{equation}\label{remint}
\frac{2}{\pi} d \int_0^{k_0} \chi(k) R(k) \Big( \Id - \Pi_0\Big) \, dk.
\end{equation}
for some small $k_0>0$ where   $\chi(k)$ is a smooth function with value $1$ near $k=0$ and supported on $[0, k_0]$ (recall that in Part I we have first decomposed the integral near and away from $k=0$, then analyzed the part $[k_0,\infty)$ 
and dealt with the part $R(k)\Pi_{<}$ where $\Pi_{<}$ is the projector onto the part correponding to the negative spectrum).

Rather than study (\ref{remint}) directly, we subtract off the free resolvent kernel, and consider
\begin{equation}
\frac{2}{\pi} d \int_0^{k_0} \chi(k) \Big( R(k) \big( \Id - \Pi_0 \big) - \phi R_0(k) \phi \Big) \, dk.
\label{r-diff}\end{equation}

First consider the $\Pi_0$ term, $\chi(k) R(k) \Pi_0 = \chi(k) k^{-2} \Pi_0$. This term does \emph{not} lie in the calculus since it is not rapidly vanishing at $\lb$, $\rb$ or $\bfc$. To deal with this we introduce 
the function $\chi((x+x')k/xx')$, where $\chi$ is as above,  which is supported away from these boundary faces. Then $k^{-2} \chi((x+x')k/xx') \Pi_0$ lies in the calculus. On the other hand, the remaining part $k^{-2} (1 - \chi((x+x')k/xx') \Pi_0$ can be treated directly: we have
$$
\frac{2}{\pi} d \int_0^{k_0} \chi(k) k^{-2} (1 - \chi((x+x')k/xx')  \, dk \leq C \frac{x+x'}{xx'}.
$$
So consider the kernel $(x+x')/xx' \Pi_0$ on $M^2$. Let $s = x/x'$. With respect to the scattering (Riemannian) half-density $|dg dg'|^{1/2}$ on $M^2$, $\Pi_0$ decays as $(xx')^{n-2+m'}$. We can write the kernel $d \big( (x+x')/xx' \Pi_0 \big)$ as $A_1 + A_2$ where $A_1$ is supported in $s \geq 1$ and $A_2$ is supported in $s \leq 1$. Notice that acting with $d$ on the left factor yields an extra decay factor of $x$. So $A_1$ is bounded pointwise by ${x'}^{2n-4+2m'} s^{n-1+m'} = {x'}^{n+\delta} s^{\alpha + \delta}$ with $\delta = n-4+2m'$, $\alpha = 3-m'$. Notice that $\delta > 0$ by our assumptions on $m'$ relative to $n$. The argument in the proof of Proposition 5.1 of Part I then shows that $A_1$ is bounded on $L^p$ for $p < n/\alpha$. Similarly, $A_2$ has kernel bounded by $x^{2n-4+2m'} s^{-(n-2+m')} = x^{n+\delta} s^{-(\beta + \delta)}$, $\beta = 2-m'$ and is bounded on $L^p$ for $p > n/(n-\beta)$. Therefore the kernel $A_1 + A_2$ satisfies the conditions of the theorem.

We return to \eqref{r-diff}, replacing $R(k) \Pi_0$ with $k^{-2} \chi((x+x')k/xx') \Pi_0$ as allowed by the discussion above. 
By Theorem~\ref{res-conic-1} (or Theorem~\ref{thm:nullspace} in the asymptotically Euclidean case), $R(k)$ is in the calculus with pseudodifferential order $-2$ and with index sets $-2$ at $\zf$, $n/2 - 4 + m'$ at $\rbo$ and $\lbo$, $-2$ at $\bfo$, $0$ at $\sca$ and trivial at all other boundary hypersurfaces. Also, we have just seen that $\chi((x+x')k/xx')k^{-2} \Pi_0$ is in the calculus; it has pseudodifferential order $-\infty$ and index sets $-2$ at $\zf$, $n/2 - 4 + m'$ at $\rbo$ and $\lbo$, $n-6+2m'$ at $\bfo$ and trivial at all other boundary hypersurfaces. Of course, the $\Pi_0$ term is chosen precisely to cancel the leading order behaviour of $R(k)$ at $\zf$, and by  Theorem~\ref{res-conic-1}, the leading behaviour of the difference 
$R(k) - \chi((x+x')k/xx')k^{-2} \Pi_0$ is at order $> -1$. 
Notice also that, with our assumptions on $m'$ relative to $n$, $n-6+2m' > -1$ and so the leading behaviour of $R(k) - \chi((x+x')k/xx')k^{-2} \Pi_0$ is equal to that of $R(k)$, at order $-2$,  with the next term at order $-1$. 

Now we bring in the kernel $\phi R_0(k) \phi $. This is in the calculus with pseudodifferential order $-2$ and with similar index sets (actually, better at $\lbo$ and $\rbo$ if $m < 2$). The point of subtracting $\phi R_0(k) \phi $ is that it \emph{cancels the leading behaviour of $R(k) - \chi((x+x')k/xx')k^{-2} \Pi_0$ at $\bfo$ and at $\sca$}, so that 
$R(k) - \chi((x+x')k/xx')k^{-2} \Pi_0 - \phi R_0(k) \phi $ has leading behaviour at order $\geq -1$ at $\bfo$ and $\geq 1$ at $\sca$.

Now we apply $d = \sum_i dz_i \otimes  \partial_{z_i}$ on the left. Vector fields $\partial_{z_i}$ of finite length on $M$ lift from the left factor to have the form $\rho_{\bfo}\rho_{\lbo}\rho_{\lb}\rho_{\bfa}$ times a vector field tangent to the boundary of $\MMksc$. Therefore 
$d$ increases the order of vanishing to $n/2 - 3 + m'$ at $\lbo$ and to $0$ at $\bfo$ and leaves the others fixed (and also increases the pseudodifferential order to $-1$). 
We can now apply Proposition 5.1 of \cite{GH}, with $\alpha = 3-m'$ and $\beta = 2-m'$,  to deduce that \eqref{r-diff} is bounded for the stated range. Combined with the earlier results about the integral for large $k$, we have proved that $T - \phi T_0 \phi$ is bounded on $L^p$ for the stated range, where $T_0$ is the classical Riesz transform on the exact cone (with respect to $\Delta_0$). However, it is shown in \cite{Li} that $T_0$ is bounded on $L^p$ for $1 < p < \infty$ (this is of course classical in the asymptotically Euclidean case where $T_0$ is the Riesz transform on $\RR^n$), so we conclude that $T$ itself is bounded on $L^p$ for the stated range. \\

Finally we prove that the range of $p$ is sharp. For simplicity we shall do this first in the case that $m'<2$. In this case we can write the resolvent kernel near the right boundary (now as a multiple of the scattering half-density rather than the b-half density, which gives an extra factor of $(xx')^{n/2}$) as 
\begin{equation*}\begin{gathered}
R(k) = {x'}^{\ndemi} k^{\ndemi -4+m'}\frac{{\kappa'} K_{\ndemi-1+m'}(\kappa')}{\Gamma(\ndemi-1+m')2^{\ndemi+m'-2}}x^{n/2-1}\sum_{j=0}^Nb^j_{\nu_m}(y')\varphi_j(z) \\ 
= {x'}^{n-4+m'} {\kappa'}^{\ndemi -4+m'}\frac{{\kappa'} K_{\ndemi-1+m'}(\kappa')}{\Gamma(\ndemi-1+m')2^{\ndemi+m'-2}}x^{n/2-1}\sum_{j=0}^Nb^j_{\nu_m}(y')\varphi_j(z) 
\end{gathered}\end{equation*}
modulo $O(\rho_{\rbo}^{n - 3 + m'} \log\rho_{\rbo} )$, with $b_{\nu_m}$ some non-zero function on $\pl M$
in $E_{\nu_m}$.  So we have 
\begin{multline*}
d R(k) - d k^{-2} \Pi_0  = {x'}^{n-4+m'}\Big( \frac{{\kappa'}^{\ndemi - 3 + m'} K_{\ndemi-1+m'}(\kappa')}{\Gamma(\ndemi-1+m')2^{\ndemi+m'-2}} - \kappa^{-2} \Big) \\ \times \sum_{j=0}^Nb^j_{\nu_m}(y')d (x^{n/2-1} \varphi_j(z)) + O({x'}^{n-3+m'} \log x')
\end{multline*}
(since $d$ acts in the left variable $z$ here). Integrating in $k$, we change variable to $\kappa' = k/x'$ (which gives us an extra power of $x'$) and obtain 
$$
 {x'}^{n-3+m'} \int_0^\infty  \Big( \frac{{\kappa'}^{\ndemi - 3 + m'} K_{\ndemi-1+m'}(\kappa')}{\Gamma(\ndemi-1+m')2^{\ndemi+m'-2}} - \kappa^{-2} \Big) \, d\kappa'  \sum_{j=0}^N \frac{b^j_{\nu_m}(y')d (x^{n/2-1} \varphi_j(z))}{ \Gamma(\ndemi-1+m')2^{\ndemi+m'-2}}
$$
for the leading asymptotic of $T$ at $x' = 0$. 
Since $x^{n/2 - 1} \varphi_j(z)$ is not constant, $d (x^{n/2 - 1} \varphi_j(z))$ is not identically zero. Moreover, by Lemma~\ref{comp}, the $\kappa'$ integral does not vanish. Hence $T = a(z, y') {x'}^{n-3+m'} + O({x'}^{n-2+m'})$ at $x' = 0$ where $a$ does not vanish. It follows immediately that the upper threshold for $p$ is sharp. A similar analysis at the left boundary shows that $T = b(y, z') x^{n-2+m'} + O(x^{n-1+m'})$ as $x \to 0$ where $b$ does not vanish, showing that the lower threshold is also sharp. 

For $m' = 2$ we have only to verify that the extra term in \eqref{leadingrbo} does not cancel the first term after integration in $\kappa'$. But this is straightforward since the two terms have different asymptotics as $x \to 0$. 

This completes the proof of the theorem.
\qed\\

\textsl{Remark}: it is interesting to see that, for instance in the Euclidean setting and dimension $n\geq 4$ with $m<2$, if we  replace $d$ in the definition 
of Riesz transform by the first order differential operator
\[D=x(x\pl_x-(n-3+m)+Q(\pl_y)) \] where $Q(\pl_y)$ is a first order differential operator on $\pl M=S^n$ 
satisfying 
\[Q(\pl_y)\phi=0, \forall \phi\in \ker(\Delta_{S^{n-1}}-m(m+n-2)),\] 
the operator $D\circ(P_>)^{-\demi}$
is bounded on $L^p(M)$ for $n/(n-1+m)<p<n/(3-m)$.

\end{document}

%% file: riesz2.pstex_t
\begin{picture}(0,0)%
\includegraphics{riesz2.pstex}%
\end{picture}%
\setlength{\unitlength}{3947sp}%
\begingroup\makeatletter\ifx\SetFigFont\undefined%
\gdef\SetFigFont#1#2#3#4#5{%
  \reset@font\fontsize{#1}{#2pt}%
  \fontfamily{#3}\fontseries{#4}\fontshape{#5}%
  \selectfont}%
\fi\endgroup%
\begin{picture}(6801,6611)(2834,-7355)
\put(4643,-4666){\makebox(0,0)[lb]{\smash{{\SetFigFont{12}{14.4}{\rmdefault}{\mddefault}{\updefault}{$\sca$}%
}}}}
\put(4573,-3056){\makebox(0,0)[lb]{\smash{{\SetFigFont{12}{14.4}{\rmdefault}{\mddefault}{\updefault}{$\lb$}%
}}}}
\put(6323,-5806){\makebox(0,0)[lb]{\smash{{\SetFigFont{12}{14.4}{\rmdefault}{\mddefault}{\updefault}{$\rb$}%
}}}}
\put(6513,-3246){\makebox(0,0)[lb]{\smash{{\SetFigFont{12}{14.4}{\rmdefault}{\mddefault}{\updefault}{$\bfo$}%
}}}}
\put(8453,-3116){\makebox(0,0)[lb]{\smash{{\SetFigFont{12}{14.4}{\rmdefault}{\mddefault}{\updefault}{$\zf$}%
}}}}
\put(8303,-4566){\makebox(0,0)[lb]{\smash{{\SetFigFont{12}{14.4}{\rmdefault}{\mddefault}{\updefault}{$\rbo$}%
}}}}
\put(6753,-1866){\makebox(0,0)[lb]{\smash{{\SetFigFont{12}{14.4}{\rmdefault}{\mddefault}{\updefault}{$\lbo$}%
}}}}
\put(4183,-5466){\makebox(0,0)[lb]{\smash{{\SetFigFont{12}{14.4}{\rmdefault}{\mddefault}{\updefault}{$\bfa$}%
}}}}
\end{picture}%